\documentclass[runningheads]{llncs}
\usepackage{graphicx}
\usepackage{color}
\usepackage{amsmath, amssymb}
\usepackage{algorithmic, algorithm}
\usepackage{multirow}
\def \RR {{\mathbb R}}

\newtheorem{Proposition}{Proposition}

\begin{document}
\title{Constrained reachability problems for a planar manipulator}
\author{Simone Cacace\inst{1}\and
Anna Chiara Lai\inst{2}\orcidID{0000-0003-2096-6753} \and
Paola Loreti\inst{2}}
\authorrunning{S. Cacace et al.}
\institute{Dipartimento di Matematica e Fisica, Universit\`{a} degli studi Roma Tre, Rome, Italy \\
\email{cacace@mat.uniroma3.it} \and
Dipartimento di Scienze di Base e Applicate per l'Ingegneria, Sapienza Universit\`{a} di Roma, Rome, Italy\\
\email{\{anna.lai,paola.loreti\}@sbai.uniroma1.it}}
\maketitle              % typeset the header of the contribution
\begin{abstract}
We address an optimal reachability problem for a planar manipulator in a constrained environment. 
After introducing the optmization problem in full generality, we practically embed the geometry of the workspace in the problem, by considering some classes of obstacles. To this end, we  present an analytical approximation of the distance function from the ellipse.   
We then apply our method to particular models of hyper-redundant and soft manipulators, by also presenting some numerical experiments.  

\end{abstract}
\section{Introduction}
We address an optimal reachability problem for a planar manipulator in a constrained environment, 
%We analyze three fundamental tasks: to avoid (a possibly disconnected) compact obstacle in the working space, to reach a target point in the plane with the  end-effector, and to minimize a  quadratic cost on the actuated controls. 
 which is part of an ongoing investigation on snake-like and octopus-like manipulators in the framework of optimal control theory of partial differential equations. The  models discussed in the present paper were originally introduced in \cite{CLL18}. Subsequent works by the authors refined the investigations in the cases of uncontrolled regions of the manipulators (modeling mechanical breakdowns)\cite{CLL19}, and grasping tasks \cite{CLL19ext}. Part of the results presented here earlier appeared in \cite{CLL20}, the main novelty consists in a more general setting of the problem, the investigation of a much wider class of  obstacle geometries, and related new numerical tests.

Our setting is stationary, namely we optimize the shape of a planar manipulator at the equilibrium. 
 We begin our investigation by considering the problem in full generality, from an optimal control theoretic perspective. We introduce a cost functional encompassing by penalization the obstacle avoidance and reachability tasks, and a quadratic running cost on the controls. The problem is then to minimize such functional in the set all of the admissible equilibrium configurations of the manipulator.  
 
 Then we address the issue of practically encompassing the geometry of the working space in the problem. More precisely, the obstacle avoidance task is enforced by introducing an elastic potential steering the manipulator outside the obstacles. In our setting, such potential is deeply related to the distance function from the obstacles. Our case study includes obstacles composed by circles, squares and ellipses. In particular, the study of the distance function from an ellipse involves root finding of a quartic polynomial, and its numerical  computation can result cumbersome in the case of general or time-varying ellipses \cite{p2e}. We present an analytical approximation of the distance function from the ellipse,  based on the linearization of an explicit formula for the roots of quartic polynomials. Moreover, we describe the approximation of the distance function for general closed obstacles with compact boundary.
 
 In the second part of the paper, we specialize the optimization problem to the case of two classes of planar manipulators: a hyper-redundant manipulator and a soft manipulator. These devices share the same physical features, respectively declined in an either discrete or continuous fashion.  We assume indeed an inextensibility constraint, a non-uniform  angle/curvature  constraint, a  bending moment (on the joints in the discrete case and pointwise in the continuous one) and angle/curvature controls. A Lagrangian formulation of the dynamics is introduced for both models, and we provide an explicit characterization of the equilibria. Finally, the optimal reachability problem with obstacle avoidance  is numerically solved in some cases of interest. 
 
From the seminal paper \cite{hyper}, where the hyper-redundant manipulators were firstly introduced, countless papers were devoted to the control of octopus-like manipulators in constrained environments, see for instance \cite{multi,nature,laschi1,laschisurvey} and the reference therein for a general introduction. The papers that mostly inspired our work include \cite{chiri}, for an early study on the interplay between the continuous and discrete settings, and \cite{optimal,rev4} for an optimal control theoretic approach to constrained reachability problems. We also refer to the papers  \cite{kinematicssoft,dynamicoctopusrobot,hand,fibonacci,softinspiredoctopus} for a modeling overview.  
 
\paragraph{Organization of the paper.}
In Section \ref{s2}, we introduce the optimal control problem, while Section \ref{s3} is devoted to the computation of  distance functions from compact sets. In Section \ref{s4} and Section \ref{s5}, we specialize the optimal control problem, respectively to a class of hyper-redundant and soft manipulators, and we present some numerical simulations.

\section{A general optimal control problem for constrainded reachability}\label{s2}
In this section, we consider a general, unidimensional planar manipulator, whose stationary configuration is  modeled by a function $q(s;u):=[0,1]\to\mathbb R^2$ depending on its arclength coordinate $s$ and on a \emph{control}  $u:A\subseteq[0,1]\to U$, where $U$ is the \emph{control set}. The function $q$ is described as a solution of an either controlled continuous differential equation or a difference equation, in the form 
\begin{equation}\label{diffeq}
\begin{cases} q^\prime=f(q,u)\\
q(0)=q_0\in\RR^2\,,\end{cases}\end{equation}
where, with a little abuse of notation, $q^\prime$ denotes either a derivative or a finite difference, and $f:\RR^2\times U\to \RR^2$. 
However, in the special cases treated in the present paper,  we also have an explicit input-to-state map $u\mapsto q(\cdot;u)$. The domain $A$  of the control function depends on the adopted model. For instance, if we are dealing with a discrete manipulator, $A$ is a finite (or countable) set of points corresponding to the joints. Otherwise, if we are dealing with a soft robot, we may set $A=[0,1]$, meaning that the controls are enforced pointwise on the whole manipulator. The set $A$ may be a finite union of intervals to model scenarios in which only a portion of the manipulator is controlled, see for instance \cite{CLL19}. We denote by $\mathcal A$ the set of \emph{admissible configuration-control pairs}, that is the couples $(q,u)$ such that $u$ is a control function, $q(s)=q(s;u)$ is the corresponding configuration, and such that some regularity assumptions are satisfied. For instance, if we are in a continuous, differential setting, one can define $\mathcal A$ as 
{$$\mathcal A:=\{ (q,u)\mid \text{\small $u:[0,1]\to U$ is measurable and $q$ is a Carath\'eodory solution of \eqref{diffeq}}\}.$$}

Concerning the working space geometry, we denote by $\Omega\subset\RR^2$ a closed subset of $\RR^2$ with compact boundary representing an \emph{obstacle}. In our examples $\Omega$ is either a circle, a square, an ellipse or a finite union of these objects. We take into account also the \emph{distance function} from $\Omega^c:=\RR^2\setminus \Omega$ 
$$q\mapsto \mathbf d(q,\Omega^c):=\inf_{x\in \Omega^c}\{|x-q|\}.
 $$
 The \emph{target} is a point $q^*\in \RR^2$. Finally, we consider a \emph{running cost} $\ell(q,u):\RR^2\times U\to [0,+\infty)$. For instance, a \emph{quadratic cost on the controls} is independent from the position of the manipulator, and it reads $\ell(u)=u^2$.
 
In this setting, we consider the problem of finding an admissible configuration-control pair $(q,u)$ such that
\begin{enumerate}
%    \item $q(s)=q(s;u)$ i.e., $q$ is the configuration corresponding to the control $u$;
    \item $q$ avoids the obstacle $\Omega$ minimizing the tip-target distance $|q(1)-q^*|$;
    \item $(q,u)$ minimizes the associated \emph{integral cost} $$\int_0^1\ell(q(s),u(s))ds.$$ 
\end{enumerate} 
The problem can be attacked by considering the cost functional:
 \begin{equation}\label{contactfunctional}
\begin{split}
    \mathcal  J(q,u):=&\frac12\int_0^1 \ell(q(s),u(s)) ds+\frac{1}{2\delta}|q(1)-q^*|^2\\&+ \frac{1}{2\tau}\int_0^1 
 \mathbf d^2(q(s),\Omega^c)ds,
\end{split}\end{equation}
with penalty parameters $\delta,\tau>0$.
We recognize in the first two terms of $\mathcal J$ the integral cost and the tip-target distance. The third term vanishes if and only if there is no interpenetration of $q$ with the obstacle $\Omega$, i.e.,  this term encompasses the obstacle avoidance task as $\tau\to0$. 
Then, we recast the optimal reachability problem as the following constrained optimization problem
\begin{equation}\label{pbgen}
\text{minimize } \mathcal J(q,u) \quad \text{ subject to } (q,u)\in \mathcal A.
\end{equation}
In Section \ref{s4} and \ref{s5}, we specialize this problem to a class of hyper-redundant and soft manipulators, by providing an explicit description of the underlying control model. 

 \section{Exact and approximated distance formulas  for obstacles}\label{s3}
  In this section, we collect some distance formulas for obstacles with compact boundary. In our tests we take into exam circular, square and elliptic obstacles. We recall here the distance functions of a point $q=(q_1,q_2)$ from the boundaries of a square of side $l$ and of a circle of radius $r$, centered in $c\in \RR^2$:
 $$\mathbf d_{square}^2(q):= \min_{a_1,a_2\in\{0,-l/2,l/2\}} \min\{|q-c+(a_1,a_2)|^2\}$$
 $$ \mathbf d_{circle}^2(q):=(|q-c|-r)^2.$$ 
 In what follows, we take into exam an analytical approximation for the distance function from an ellipse, and we describe a strategy for the numerical approximation of the distance function from general sets with compact boundaries. 
 \subsection{Distance formulas from the ellipse}
 
 Let $0<b\leq a$ and consider the ellipse $\mathcal E(a,b)$ centered in the origin with semi-axes $a$ and $b$, implicitly defined by the equation 
\begin{equation}\label{Edef}
 E(x):=E(x_1,x_2)=\left(\frac{x_1}{a}\right)^2+\left(\frac{x_2}{b}\right)^2-1=0.\end{equation}
 We define the square distance from $\mathcal E(a,b)$ by
$$\mathbf d^2(q):=\min\{|x-q|^2\mid E(x)=0, x\in\RR^2\}.$$
We fix $q\in\RR^2$ and we  use the Lagrange multiplier method to investigate $\mathbf d(q)$. Consider the Lagrangian function
$$L_q(x,\lambda):=|x-q|^2-\lambda E(x).$$
The minimization of $L_q$ leads to the optimality system with unkowns $x=(x_1,x_2)$ and $\lambda$:
\begin{equation}
\begin{cases}q_1=x_1\left(1-\frac{\lambda}{a^2}\right)\\
q_2=x_2\left(1-\frac{\lambda}{b^2}\right)\\
E(x)=0\\
\lambda\leq b^2.\end{cases}
\end{equation}
Note that the first three equations are stationarity conditions, while the inequality in the multiplier $\lambda$ is an actual local minimality condition. By algebraic computations one ends up with the equivalent formulation:
\begin{equation}\begin{cases}
q_1=x_1\left(1-\frac{\lambda}{a^2}\right)\\
q_2=x_2\left(1-\frac{\lambda}{b^2}\right)\\
P(\lambda):=\left(\left((\lambda-a^2\right)\left(\lambda-b^2\right)\right)^2-a^2q_1^2\left(\lambda-b^2 \right)^2-b^2q_2^2\left(\lambda-a^2\right)^2=0 \\
\lambda\leq b^2.\end{cases}
\end{equation}
Now, one can prove that the required multiplier $\lambda^*(q)$ is the smallest root of $P$. Indeed, the case $a=b$ corresponds to the circle, and it is trivial to check that $\lambda^*(q)=a^2-a|q|\leq a^2=b^2$. If otherwise $b<a$ and if $q_2\not=0$,  
 %\in R_{a,b}:=\{(q_1,0)\mid |q_1|< \frac{a^2-b^2}{a}\}$ 
 then $\lambda^*(q)$ is univoquely determined by the above system, since $P$ admits one and only one root in the interval $(-\infty,b^2)$. Finally, if $b<a$ and $q_2=0$, then $P(b^2)=P(a^2-a|q_1|)=0$ and a direct computation implies the global minimum of $L_q$ to be attained at points of the form $(x^*(q),\lambda^*(q))$ with $\lambda^*(q)=\min\{a^2-a|q_1|, b^2\}$. 
 Hence, the exact formula for the distance is given by
 \begin{equation}\label{dexact}\mathbf d^2(q)=\begin{cases}
 \left(\dfrac{\lambda^*(q)}{a^2-\lambda^*(q)}q_1\right)^2+ \left(\dfrac{\lambda^*(q)}{b^2-\lambda^*(q)}q_2\right)^2&\text{ if } q_2\not=0\\& \text{ or } a^2-a|q_1|< b^2\\\\
 b^2&\text{ if } q_1=q_2=0\\\\
  b^2-\dfrac{b^2}{a^2-b^2}q_1^2&\text{ if } q_2=0, q_1\not=0\\& \text{ and } a^2-a|q_1|\geq b^2.
 \end{cases}\end{equation}
 Note that, if $b<a$, the second case in above expression is a particular case of the third case. 
 Our idea is to use explict formula for the roots of quartic polynomials to approximate $\lambda^*$ as $\varepsilon:=a^2-b^2\to 0^+$. Let us rewrite $P(\lambda)$ as 
 $$P_\varepsilon(\lambda):=\left(\left((\lambda-b^2-\varepsilon\right)\left(\lambda-b^2\right)\right)^2-(b^2+\varepsilon)q_1^2\left(\lambda-b^2 \right)^2-b^2q_2^2\left(\lambda-b^2-\varepsilon\right)^2\,.\\
 $$
Clearly, $\lambda^*(q)$ is also the smallest root of $P_\varepsilon$, and we denote by $\lambda_\varepsilon(q)$ its first order approximation, so that $\lambda^*(q)=\lambda_\varepsilon(q)+o(\varepsilon)$ as $\varepsilon\to 0^+$ for all $q\in \RR^2$. Then we replace $\lambda^*(q)$ in \eqref{dexact} by $\lambda_\varepsilon(q)$:
 \begin{equation}\label{dplug}\bar{\mathbf  d}_\varepsilon ^2(q):=\begin{cases}
 \left(\dfrac{\lambda_\varepsilon(q)}{b^2+\varepsilon-\lambda_\varepsilon(q)}q_1\right)^2+ \left(\dfrac{\lambda_\varepsilon(q)}{b^2-\lambda_\varepsilon(q)}q_2\right)^2&\text{ if } q_2\not=0 \text{ or } q_1^2> \dfrac{\varepsilon^2}{b^2+\varepsilon}\\
% b^2&\text{ if } q_1=q_2=0\\
  b^2-\dfrac{b^2}{\varepsilon}q_1^2&\text{ otherwise. }
 \end{cases}\end{equation}
 Incidentally, notice that $\bar{\mathbf d}_\varepsilon^2(q)\leq b^2-\varepsilon \frac{b^2}{b^2+\varepsilon}$ when $q_2=0$ and $q_1^2\leq \dfrac{\varepsilon^2}{b^2+\varepsilon}$.
We performed a symbolic computation using the Wolfram Mathematica software to get the following first order approximation of $\bar {\mathbf d}^2_\varepsilon$:
\begin{equation}\label{dapprox}{\mathbf  d}_\varepsilon ^2(q):=\begin{cases}
 (b-|q|)^2+\varepsilon \frac{q_1^2}{|q|}(b-|q|)&\text{ if } q_2\not=0 \\&\text{ or } q_1^2\geq \dfrac{\varepsilon^2}{b^2+\varepsilon}
  \\\\
    b^2-\dfrac{b^2}{\varepsilon}q_1^2&\text{ otherwise. }
\end{cases}\end{equation}
 By construction, we finally get, for all $q\in \RR^2$,  the estimate
 $$\mathbf d^2(q)=\mathbf d^2_\varepsilon(q)+o(\varepsilon) \quad \text{ as } \varepsilon\to 0^+.$$
% \paragraph{General ellipses} A generic ellipse $\mathcal E(A,B)$ is implicitely defined by the equation
%$$E_{A,B}(x):= x^t A x+B^t x-1=0$$
%where $B\in \RR^2$, $A\in \RR^2\times \RR^2$, $A$ is symmetric and $det(A)>0$. General ellipses can be reduced to canonical forms by a change of variable. In particular since $A$ is symmetric, then there exist an orthogonal matrix $M$, a vector $C$ and two real numbers $0<b\leq a$ such that every $\bar x\in \mathcal E(A,B)$ satisfies 
%$$\bar x = M x+C$$ 
%for some $x\in \mathcal E(a,b)$, see \eqref{Edef}. Therefore the distance function $\bar {\mathbf d}$ from $\mathcal E(A,B)$ satisfies
%$$\bar {\mathbf d}(q)=\mathbf d(M^t(q-C)) $$
\subsection{Distance function from general obstacles}
\begin{figure}[ht!]
	\centering
	\begin{tabular}{c}
		\includegraphics[width=\textwidth]{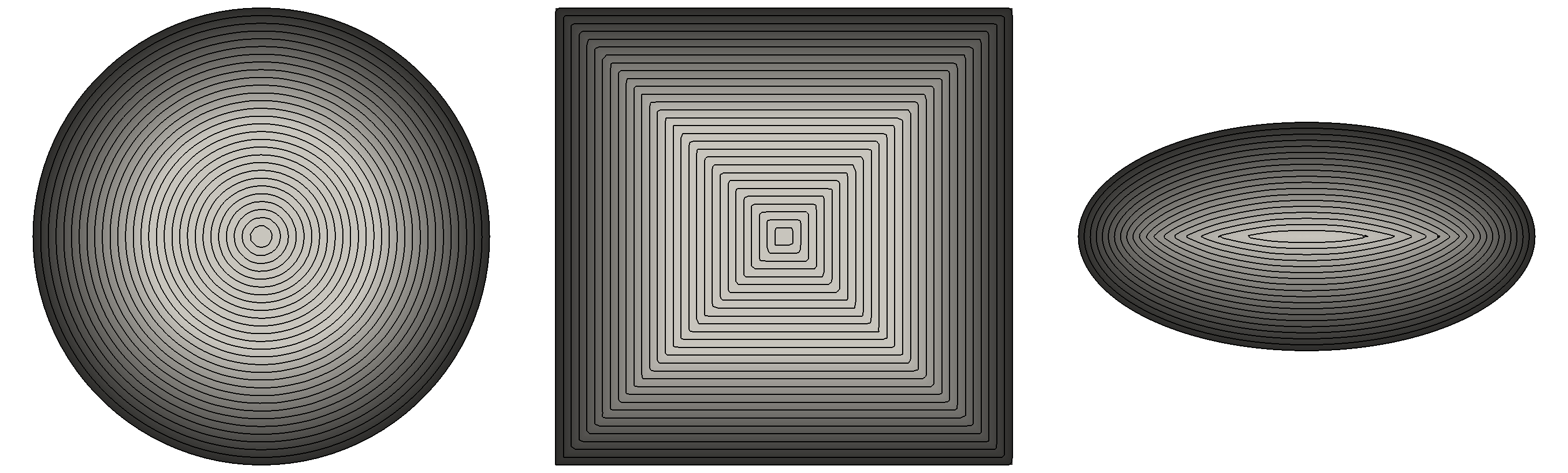}
	\end{tabular}
		\caption{\label{distances}  Level sets of the distance function for a circle, a square, an ellipse.} 
\end{figure}

When dealing with a general obstacle, analytical expressions for the distance function are no longer available. Nevertheless, from a theoretical point of view, the distance function can be characterized as the solution of a first order partial differential Hamilton-Jacobi equation, the celebrated Eikonal equation:
$$
\left\{
\begin{array}{ll}
     |\nabla \mathbf d(x)|=1 &\quad x\in \Omega \,, \\
     \mathbf  d(x)=0&\quad x\in\RR^2\setminus \Omega\,. 
\end{array}
\right.
$$
It is well known that the distance function is merely continuous, since its gradient can exhibit singularities. This is the case even for the examples discussed above, namely the distance function for the circle is not differentiable at its center, for the square on the diagonals, and for the ellipse on the segment joining its foci (see Figure \ref{distances}).
Hence, the solution to the Eikonal equation should be meant in a suitable weak sense, introducing the notion of viscosity solutions. There is a wide literature on this subject, also from a numerical point of view, which dates back to the seventies and it is still growing nowadays. This is far beyond the scope of the present paper, and we refer the interested reader to \cite{lions} as a starting point. Here, we just remark that the Eikonal equation can be solved numerically employing one of the available state-of-the-art algorithms, such as the fast marching method (see \cite{Sethian1591}). To this end,  it is enough to provide the solver a triangulation of $\Omega$, and  impose the Dirichlet condition $\mathbf d=0$ on the discrete boundary. Once the numerical solution is computed, it can be extended to the whole space via interpolation.  

\section{Optimal control of a class of hyper-redundant manipulators}\label{s4}
We consider the optimal control problem introduced in Section \ref{s2} in the case of  a planar hyper-redundant manipulator, whose joints are subject to an angular constraint, a bending moment and an angular control. This  model was earlier introduced in \cite{CLL18} and later extended to a more general setting in \cite{CLL20}. Here, after recalling the main features and properties of the model, we address the associated optimal constrained reachability problem for different types of obstacles.
%\begin{Notation}
\subsection{The model} The planar manipulator under exam is composed by $N$ rigid links and $N+1$ joints.  We denote by $m_k$ the mass of the $k$-th joint, for $k=0,\dots,N$, and we consider negligible the mass of the corresponding links. The positions of the joints are stored in the array $q=(q_0,\dots,q_N)$, where $q_0:=(0,0)$ is the anchor point. To make some of the definitions below consistent, we also consider the ghost joints $q_{-1}:=q_0+(0,\ell_0)$ for some positive $\ell_0$, and $q_{N+1}:=q_N+(q_{N}-q_{N-1})$ at the free end. 
%Moreover, given vectors $v_1, v_2\in\RR^2$, we assume standard notation for the Euclidean norm and the dot product, respectively $|v_1|$ and $v_1\cdot v_2$, 
%and we define $v_1\times v_2:=v_1\cdot v_2^\bot$, where $v_2^\bot$ denotes the clockwise orthogonal vector to $v_2$. Finally, the positive part function is denoted by $(\cdot)_+$.\\
%\end{Notation}
 The features of this manipulator are the following. 

First, we have an \emph{inextensibility constraint}, representing the fact that the links of the manipulator are rigid, therefore each couple of consecutive joints satisfies
 $ |q_{k}-q_{k-1}|=\ell_k,$ for $k=1,\dots,N$, where $\ell_k>0$ is the length of the $k$-th link. We introduce this constraint exactly, by considering the functions
\begin{equation}\label{F}
F_k(q,\sigma):=\sigma_k\left(|q_{k}-q_{k-1}|^2-\ell_k^2\right)\, \text{for } k=1,\dots,N\,,\end{equation}
where $\sigma_k$ is a Lagrange multiplier.
 
 The second matter under exam is the behavior of the joints. The model prescribes that two consecutive links, say the $k$-th and the $k+1$-th, tend to resist to bending and, however, they cannot form an angle larger in modulus than a fixed threshold $\alpha_k$. These two constraints are introduced via penalization, i.e., by considering two angular elastic potentials.
We set
 \begin{equation}\label{B}B_k(q):=\varepsilon_k b_k^2(q)\,,\end{equation}
 with
 $$b_k(q):=(q_{k+1}-q_k)\times(q_{k}-q_{k-1})\,,$$
and $v_1\times v_2:=v_1\cdot v_2^\bot$, where $v_2^\bot$ denotes the clockwise orthogonal vector to $v_2$. 
The function $B_k(q)$ represents an elastic potential, with penalty parameter $\varepsilon_k>0$, associated to the \emph{bending moment}, corresponding to the constraint $b_k(q)=0$.
Similarly, we set 
\begin{equation}\label{G}
G_k(q):=\nu_k g^2_k(q),
\end{equation}
with
 $$g_k(q):=\Bigg(\cos(\alpha_k)-\frac{1}{\ell_{k+1}\ell_k}(q_{k+1}-q_{k})\cdot(q_{k}-q_{k-1})\Bigg)_+\,,$$
 where $(\cdot)_+$ denotes the positive part of its argument. The function $G_k(q)$ is associated to the \emph{angular constraint} $g_k(q)=0$, forcing, with penalty parameter $\nu_k>0$, the relative angle between the $k$-th and $k+1$-th links in the interval $[-\alpha_k,\alpha_k]$. 
 
Finally, we consider the \emph{control term}.  We choose the control set $U:=[-1,1]$ and we prescribe the angle between the $k$-th and $k+1$-th links  to be equal to $\alpha_ku_k$ --  the control set $[-1,1]$ is chosen in order to be consistent with the angle constraint. This reduces to the  following equality constraint: 
$$b_k(q)-\ell_{k+1}\ell_k\,\sin(\alpha_k u_k)=0\,.$$
 Also in this case, we enforce the constraint via penalization, by considering
\begin{equation}\label{H}
H_k(q,u):=\mu_k\left(\ell_{k+1}\ell_k\,\sin(\alpha_k u_k)-b_k(q)\right)^2\,,\end{equation}
where $\mu_k\geq0$ is a penalty parameter. Note that to set $\mu_k=0$ corresponds to deactivate the control of the $k$-th joint and let it evolve according to the remaining constraints only. 

 We then build the Lagrangian associated to the hyper-redundant manipulator by introducing a kinetic energy term and 
{the above discussed elastic potentials:
\begin{equation}\label{LN}\begin{split}
    \mathcal L_N(q,\dot q,\sigma&,u):=\sum_{k=0}^N \frac{1}{2}m_k|\dot q_k|^2-F_k(q,\sigma)-G_k(q) -\frac{1}{2}B_k(q)-\frac{1}{2}H_k(q,u).
\end{split}\end{equation}}
%where $F_k$, $G_k$, $B_k$ and $H_k$ are respectively defined in \eqref{F},\eqref{G},\eqref{B} and \eqref{H}.

For every fixed control array $u\in [-1,1]^N$, the associated equilibria correspond to the (unique) solution of the following stationary system:  
\begin{equation}\label{eqc}
    \begin{cases}
    \nabla_{q}\mathcal L_N=0\\
        |q_{k}-q_{k-1}|=\ell_k \qquad k=1,\dots,N\\
    q_0=(0,0)\\
    q_{-1}=q_0+\ell_0(0,1)\\
    q_{N+1}=q_N+(q_{N}-q_{N-1}).
    \end{cases}
\end{equation}
We recall from \cite{CLL20} the  explicit characterization of the solutions of the above system. 
\begin{Proposition}\label{p1}
Fix $u\in[-1,1]^{N}$, assume $\alpha_k\in[0,\pi/2]$ for $k=0,\dots,N-1$, define 
$$\bar \alpha_k:=\arcsin{\frac{\mu_k}{\varepsilon_k+\mu_k}\sin{u_k\alpha_k}}$$
and, for $k=1,\dots,N$
$$ z_k:=-i\sum_{j=1}^k\ell_j e^{i\sum_{h=0}^{j-1}\bar \alpha_h}\,. $$
Then the vector $q=(q_0,q_1,\dots,q_N)$ defined by
$$q_k=
\begin{cases}
(0,0)&\quad \text{if }k=0\\
(Re(z_k),Im(z_k))&\quad \text{if }k=1,\dots,N
\end{cases}$$
is the solution of \eqref{eqc}. 
\end{Proposition}

By  Proposition \ref{p1}, if $\alpha_k\in[0,\pi/2]$ then  the  \emph{input-to-state map}
$$u\mapsto (q_1[u],\dots,q_N[u])$$
associated to  \eqref{eqc} reads
\begin{equation}\label{i2s}
    q_k[u]:=\sum_{j=1}^k \ell_j\Big( \sin{\left(\theta_j[u]\right)},-  \cos\left(\theta_j[u]\right)\Big)
\end{equation}
where
$$\theta_j[u]:=\sum_{h=0}^{j-1}\bar \alpha_h[u]; \qquad \bar \alpha_k[u]=\arcsin{\frac{\mu_k}{\varepsilon_k+\mu_k}\sin{u_k\alpha_k}}.$$

Finally, we assume  for simplicity that the total length of the manipulator is normalized to $1$, i.e., 
$\sum_{k=1}^N\ell_k=1$. Since the manipulator is composed by a series of rigid,  inextensible links, its \emph{equilibria configurations} can be parametrized by a linear interpolation $q(s;u)$ of its joints coordinates $(q_0[u],\dots,q_N[u])$:
	\begin{equation}\label{linear}q(s;u):=\frac{\ell_k- s+s_k}{\ell_k} q_{k}[u]+\frac{s- s_k }{\ell_k} q_{k+1}[u], \quad s\in (s_k, s_{k+1}],\, k=0,\dots,N-1\end{equation}
where $s_k:=\sum_{j=1}^k \ell_j$.

\subsection{Optimal reachability with obstacle avoidance}\label{optimaldisc}
We now specialize the optimal reachability problem described in Section \ref{s2} to the present model. The control set is $U=[-1,1]$, and the configuration $q(s;u)$ is given by \eqref{linear}. We choose a control quadratic running cost $\ell(u):=u^2$.  Then, for a given target point $q^*\in \RR^2$, and a closed subset $\Omega$ of $\RR^2$ representing the obstacle, problem \eqref{pbgen} reads:
\begin{equation}\label{minpb}
\min \mathcal J,\quad\text{subject to $\eqref{eqc}$ and to $u\in[-1,1]^N$\,,}
\end{equation} 
with
\begin{equation}\label{minpbJ}
    \begin{split}
    \mathcal J(q,u):=&\frac{1}{2}||u||_2^2+
    \frac{1}{2\delta}|q(1,u)-q^*|^2+\frac{1}{2\tau}\int_0^1 \mathbf d^2(q(s;u),\Omega^c)ds,   
    \end{split}
    \end{equation}
where  $||u||_2$ 
is the $l^2$ norm of the control vector $u$ and $\delta$ and $\tau$ are positive penalty parameters. Note that, due to the particular form of the input-to-state map \eqref{i2s}, the function \eqref{minpbJ} actually depends on $u$ only. 

\begin{table}[t]\centering
\begin{tabular}{|l|l|}\hline
  		Parameter description&Setting\\
  		\hline Number of links &$N=8$\\
  		Number of samples &$S=104$ ($m=13$)\\
  		\hline length of the links &$\ell_k=1/8$\\
  		bending moment &$\varepsilon_k=10^{-1}(1-0.9s_{km})$\\
  		curvature control&$\mu_k=1-0.9s_{km}$\\ 
  		penalty&\\
  		angle constraint&$\alpha_k=2\pi(2+s_{km}^2)$\\
  		\hline
  		target point&  $q^*=(0.368,-0.085)$\\
  		target penalty &$\delta=10^{-8}$\\
  		obstacle penalty &$\tau=10^{-10}$\\
  		\hline
 \end{tabular}
  	\caption{
  	Global parameter  settings for the hyper-redundant manipulator.  \label{parametertablediscr}}%}
  	\end{table}

	\subsubsection{Numerical simulations.} We discretize the parametrization interval $[0,1]$ using $S+1$ uniformly distributed samples $s_i=i/S$, for $i=0,...,S$.  Here, $S=m N$ is a multiple ($m\gg 1$) of the number of links, so that, 
	$$q(s_{km+j};u)=(1-\lambda_j)q_{k}[u]+\lambda_jq_{k+1}[u], \quad \forall k=0,\dots,N-1,\,j=0,\dots,m-1$$
	with $\lambda_j=j/m$. As in \cite{CLL20},  we approximate the integral term in \eqref{minpbJ} by a  rectangular quadrature rule, obtaining 
 a fully discrete objective function $J(u)$ with $u\in[-1,1]^N$. We then use a projected gradient descent method to solve the finite-dimensional constrained optimization of $J(u)$. Moreover, we start with $\tau\gg\delta$ and run the optimization up to convergence, then we slowly decrease $\tau$ and repeat the optimization until $\tau$ is suitably small. In this way, we first obtain an optimal configuration for the tip-target distance without considering the obstacle. Then, we iterate the procedure, to progressively penalize all the possible interpenetrations with the obstacle. In Algorithm \ref{ALG1}, we recall from \cite{CLL20} the algorithm summarizing the whole optimization process-- note that we denote by $\Pi_{[-1,1]^N}(u)$ the projection of $u$ on $[-1,1]^N$.
\begin{algorithm} 
\begin{algorithmic}[1]
\STATE{Fix $tol>0$, $tol_\tau=\tau$, and a step size $0<\gamma<1$}
\STATE{Assign an initial guess $u^{(0)}\in[-1,1]^N$}
\STATE{Compute $J(u^{(0)})$ and set $J_{tmp}=0$}
\STATE{Set $\tau>>\delta$}
\REPEAT
\STATE{$n\leftarrow 0$, $\tau \leftarrow \tau/2$}
\REPEAT
\STATE $J_{tmp}\leftarrow J(u^{(n)})$
\STATE Compute $\nabla J(u^{(n)})$
\STATE $u^{(n)}\leftarrow \Pi_{[-1,1]^N}\{u^{(n)}-\gamma \nabla J(u^{(n)})\}$
\STATE{ $n\leftarrow n+1$}
\STATE Compute $J(u^{(n)})$
\UNTIL{$|J(u^{(n)})-J_{tmp}|< tol$}
\STATE{$u^{(0)}\leftarrow u^{(n)}$}
\UNTIL{$\tau<tol_\tau$}
\end{algorithmic}
\caption{\label{ALG1}}
\end{algorithm}
The simulation parameters are summarized in Table \ref{parametertablediscr}.

\begin{table}[ht!]\centering	\begin{tabular}{|l|l|}\hline
  		Test& Obstacle\\
  		\hline
  		Test 1& $\Omega=\emptyset$\\
  		Test 2& $\Omega=B_{0.08}(0.1,-0.35)$\\
  		Test 3& $\Omega=B_{0.08}(0.1,-0.35)\cup B_{0.05}(0.3,-0.35)$\\
  		Test 4& $\Omega=Q_{0.2}^{25^\circ}(0.2,-0.35)$\\

	Test 5& $\Omega=E_{0.18,0.08}^{25^\circ}(0.2,-0.35)$\\
	Test 6& $\Omega=Q_{0.16}^{45^\circ}(0.1,-0.35)\cup E_{0.09,0.06}^{45^\circ}(0.3,-0.35)$\\
  		\hline
  	\end{tabular}
  	\caption{Obstacle settings. $B_r(x),\,Q_l^\alpha(x),\,E_{a,b}^\alpha\subset\RR^2$ denote respectively the ball of radius $r$, the square of side $l$, the ellipse of semi axes $a$, $b$, centered in $x$ and clockwise rotated by the angle $\alpha$ (in degrees). \label{controltable}}
  \end{table}
	
 \begin{figure}[ht!]
	\centering
	\begin{tabular}{c}
		\includegraphics[width=.3\textwidth]{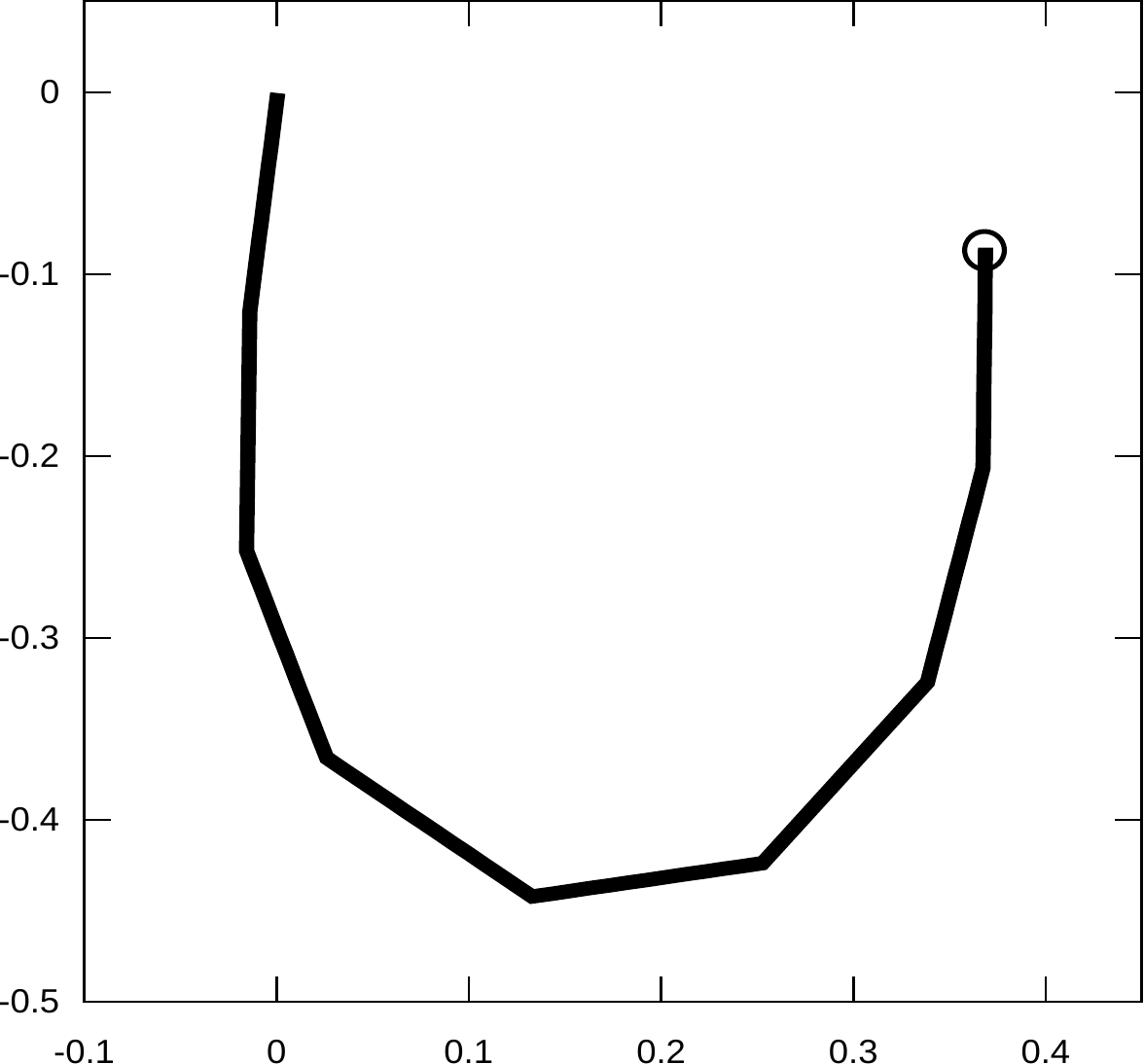}\\
	\end{tabular}
	\begin{tabular}{c}
		\includegraphics[width=.3\textwidth]{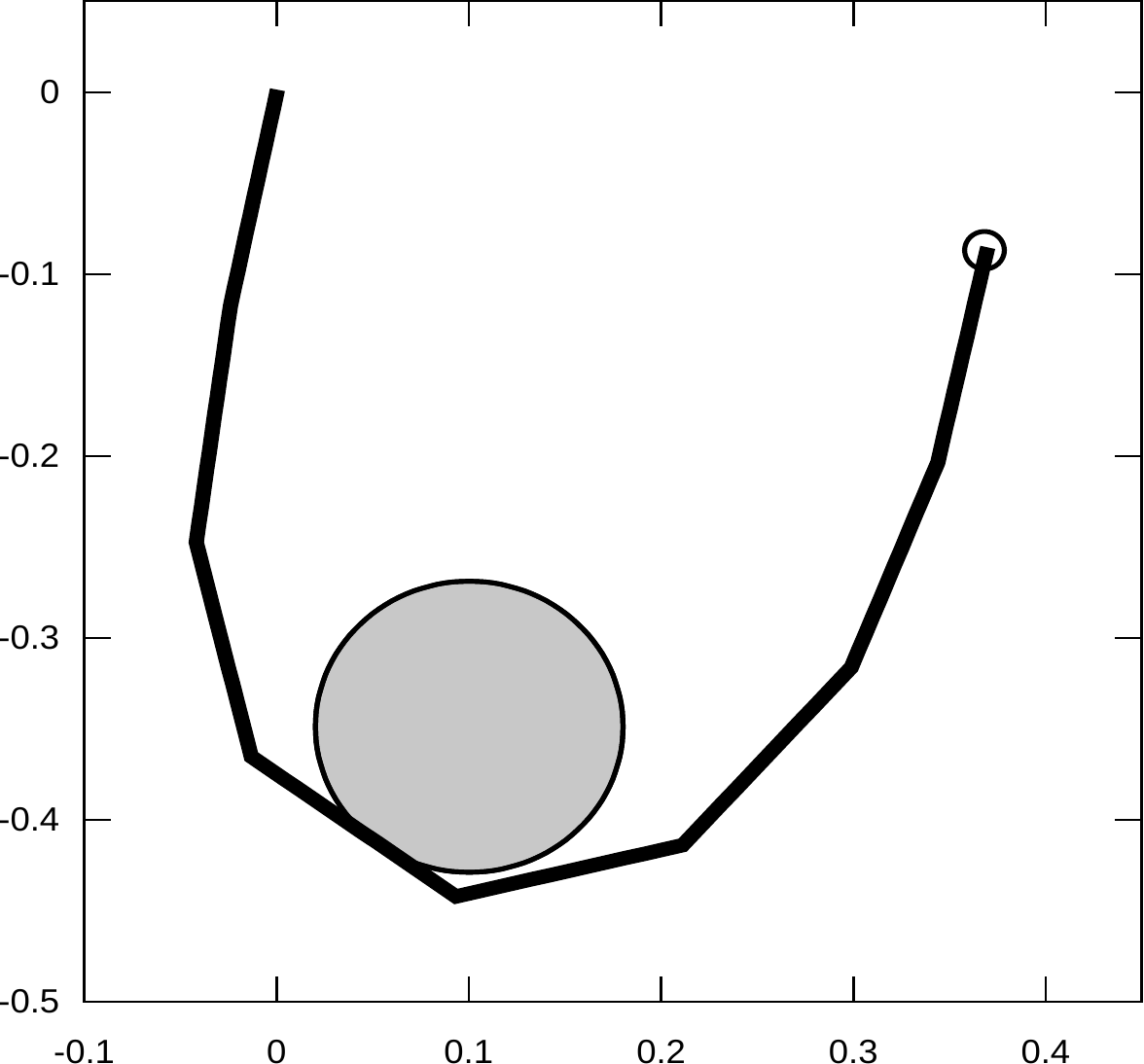}\\
	\end{tabular}
	\begin{tabular}{c}
		\includegraphics[width=.3\textwidth]{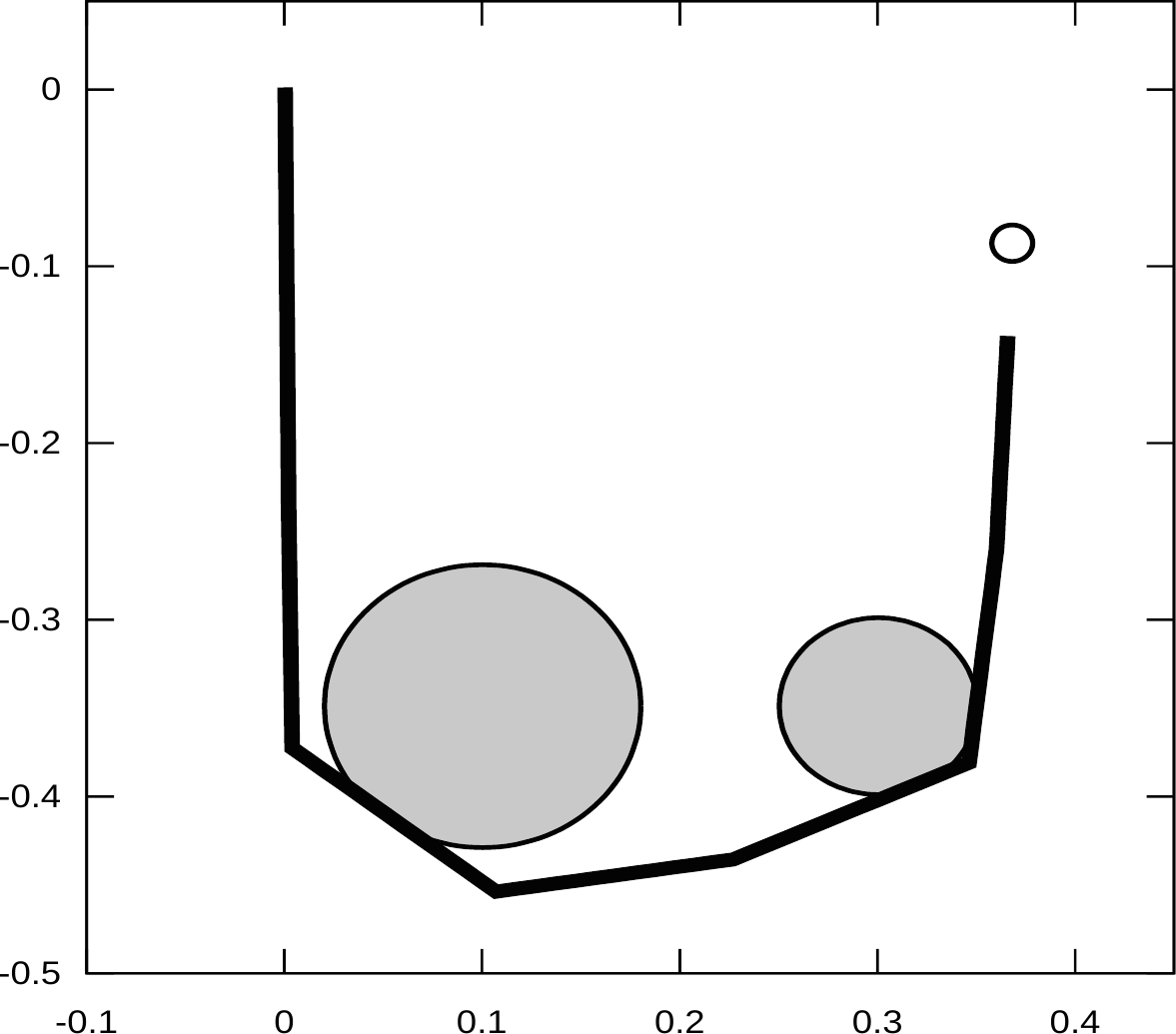}\\
	\end{tabular}
	\begin{tabular}{c}
		\includegraphics[width=.3\textwidth]{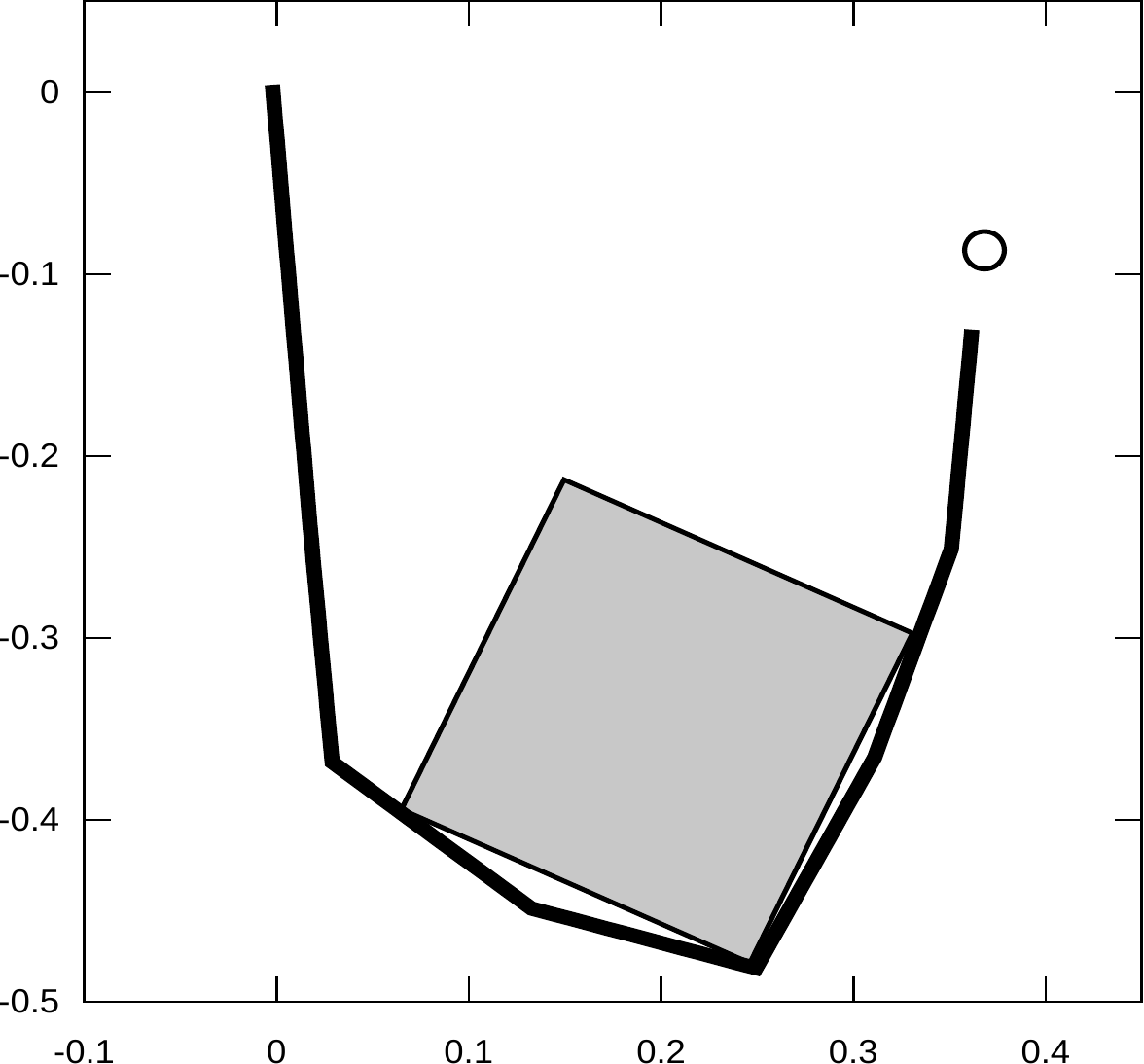}\\
	\end{tabular}
	\begin{tabular}{c}
		\includegraphics[width=.3\textwidth]{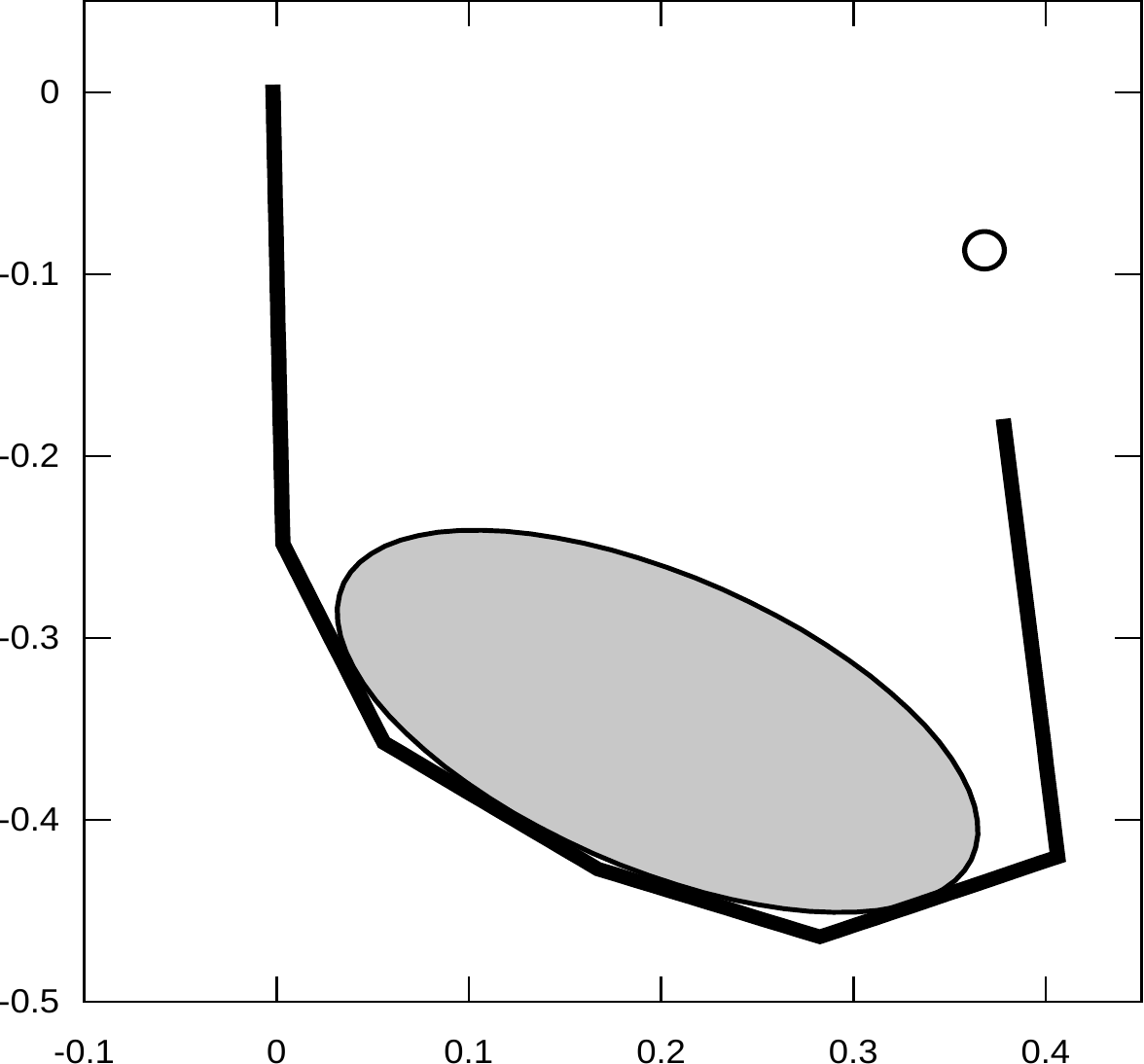}\\
	\end{tabular}
	\begin{tabular}{c}
		\includegraphics[width=.3\textwidth]{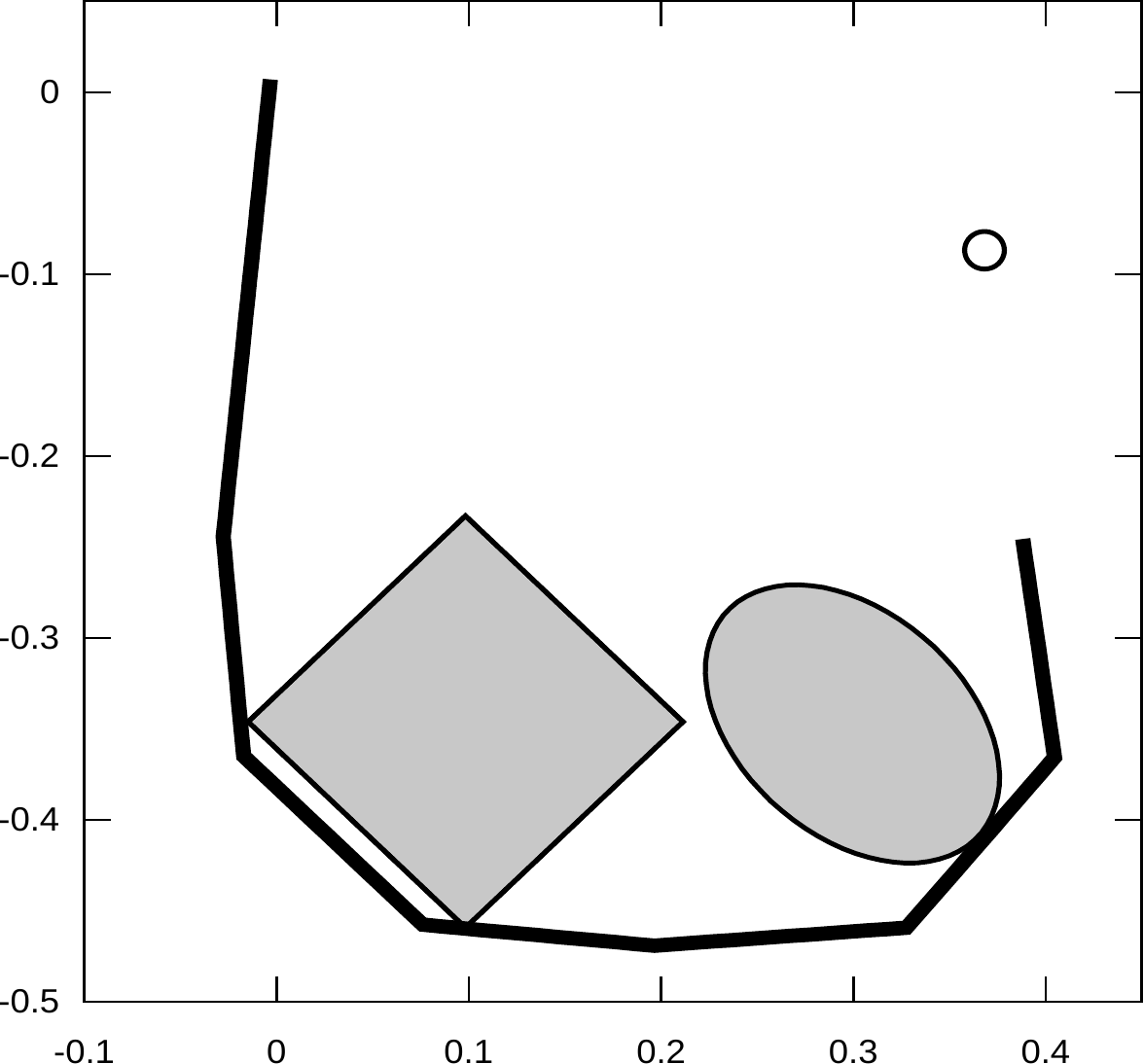}\\
	\end{tabular}
	\caption{\label{stationary3d}  The solution $q$ of Test 1-6, respectively.} 
\end{figure}

 We compare the cases reported in Table \ref{controltable}, namely the cases in which $\Omega$ is the empty set (Test 1), $\Omega$ is a ball (Test 2), $\Omega$ is the disjoint union of two balls (Test 3), $\Omega$ is a rotated square (Test 4), $\Omega$ is a rotated ellipse (Test 5), and $\Omega$ is the disjoint union of a rotated square and a rotated ellipse (Test 6).  
Note that in Test 1 and Test 2 the target $q^*$ is reached by the end-effector of the manipulator, with clearly different optimal solutions emerging from the differences between the workspaces. On the other hand, in the remaining tests, we observe that the target is unreachable, since the parameters are set in order to prioritize obstacle avoidance. 
\section{Optimal control of a class of soft manipulators}\label{s5}
 \begin{table}
\centering\caption{\label{constraints} Exact constraint equations  and associated potentials in both discrete and continuous settings.}
\scalebox{0.92}{\renewcommand\arraystretch{2}\begin{tabular}{|l|l|l|l|}\hline
 	{\bf Constraint}&&  Discrete&Continuous\\\hline
  Inextensibility&{\small Equation}& $|q_k-q_{k-1}|=\ell$&$|q_s|=1$\\\hline
\multirow{2}{*}{Curvature}&Equation&$(q_{k+1}-q_k)\cdot(q_{k}-q_{k-1})\ge \ell^2\,\cos(\alpha_k)$ &$|q_{ss}|\leq \omega$\\
&Penalization&$\nu_k\left(\cos(\alpha_k)-\frac{1}{\ell^2}(q_{k+1}-q_{k})\cdot(q_{k}-q_{k-1})\right)_+^2$& $\nu(|q_{ss}|^2-\omega^2)^2_+$\\\hline
\multirow{2}{*}{Bending}& Equation&$(q_{k+1}-q_k)\times(q_{k}-q_{k-1})=0\,$&$|q_{ss}|=0$\\
&Penalization &$\varepsilon_k\Big( (q_{k+1}-q_{k})\times(q_{k}-q_{k-1})\Big)^2$&$\varepsilon |q_{ss}|^2$\\\hline
\multirow{2}{*}{ Control}& Equation&$(q_{k+1}-q_k)\times(q_{k}-q_{k-1})=\ell^2\,\sin(\alpha_k u_k)$&$q_s \times q_{ss}=\omega u$\\
&Potential&$\mu_k\left(\sin(\alpha_ku_k)-\frac{1}{\ell^2}(q_{k+1}-q_{k})\times(q_{k}-q_{k-1})\right)^2$& $\mu\left(\omega u-q_s \times q_{ss}\right)^2$\\
\hline
  \end{tabular}}
\end{table}

We consider a soft manipulator introduced in \cite{CLL18} and encompassing the continuous counter part of the features of the hyper-redundant manipulators described in  Section \ref{s2}. In particular, the device is modeled as an inextensible elastic string subject to curvature constraints, representing a bending moment and preventing the device to bend over a fixed threshold. Moreover the curvature is forced pointwise by a control term, modeling an angular elastic internal force.
\subsection{The model} The time-varying configuration of the soft-manipulator is parametrized by the function $q:[0,1]\times [0,+\infty)\to \RR^2$.  Its evolution is determined by internal reaction forces, emerging from the inextensibility and curvature constraints and from the control term. Such  constraints and the associated angular elastic potentials are derived from the formal limit  (as the number of joints goes infinity) of the angular constraints of the hyper-redundant manipulator, see \cite{CLL19ext,CLL18}. In Table \ref{constraints}, we compare the discrete and continuous versions of the constraints under exam and the related elastic potentials. We build the continuous counter part of the Lagrangian 
introduced in \eqref{LN}:
\vskip-10pt
\begin{equation}\label{continuouslagrangian}
 \begin{split}\mathcal{L}(q,\sigma):={\int_0^1}&\Big(\underbrace{\dfrac12\rho|q_t|^2 }_\text{kinetic energy}-\underbrace{\dfrac12 \sigma(|q_s|^2 - 1) }_\text{inextensibility constr.}-\underbrace{\dfrac14\nu\left(|q_{ss}|^2-\omega^2\right)_+^2}_\text{curvature constr.}\\
 &-\underbrace{\dfrac12\varepsilon|q_{ss}|^2}_\text{bending moment}-\underbrace{\dfrac12\mu\left(\omega u-q_s\times q_{ss}\right)^2
}_\text{curvature control} \Big)ds\,,\end{split}
 \end{equation}
  where $q_t$, $q_{s}$, $q_{ss}$ denote partial derivatives in time and space respectively,  $\rho:[0,1]\to \RR^+$ is the mass distribution, $\nu,  \varepsilon, \mu:[0,1]\to\RR^+$ are the angular elastic weights associated, respectively, to the curvature constraint, the bending moment and the curvature control, while $u:[0,1]\times[0,+\infty)\to [-1,1]$ is the curvature control.

%
%
%
%
%
%In \cite{CLL18} we introduced a control model for a soft manipulator, obtained as the formal limit (as the number of the joints goes to infinity) of a sequence of the above described hyper-redundant manipulators with fixed unitary length. The mass distribution is given by the function , while the (time-dependent) configuration of the string is given by the function $q:[0,1]\times [0,+\infty)\to \RR^2$. Roughly speaking, the above introduced angular constraints and controls, prescribing the behaviuor of the relative angles of the joints, in the continuum model turn to analogous \emph{curvature} constraints and controls, see Table \ref{constraints} in which we recall from \cite{CLL19,springer19}
The equilibria of the system associated with the Lagrangian \eqref{continuouslagrangian} were explicitely characterized in \cite{CLL18}. In particular, assuming the technical condition $\mu(1)=\mu_s(1)=0$, the shape of the manipulator at the equilibrium is the solution  $q$ of
the following second order controlled ODE:
\begin{equation}\label{reducedstationary}
\quad\left\{\begin{array}{ll}
q_{ss}=\bar \omega\, u\, q_s^\bot&\mbox{in }(0,1)\\
|q_s|^2=1 &\mbox{in }(0,1)\\
q(0)=(0,0)\\q_{s}(0)=(0,-1)\,,
\end{array}
\right. 	\end{equation} 
\vskip-2pt
\noindent where  
$$ \bar\omega(s) :=\frac{\mu(s) \omega(s)}{\mu(s)+\varepsilon(s)}\,.$$ 
Assuming a sufficient regularity on the control function $u$ and solving \eqref{reducedstationary}, we obtain the following continuous version of \eqref{i2s}, namely the input-to-state map 
\begin{equation}\label{control2state}
u\mapsto q(s;u)=\int_0^s\Big(\sin(\int_0^\xi \bar\omega(z) u(z)\,dz),-\cos(\int_0^\xi \bar\omega(z) u(z)\,dz)\Big)d\xi\,.
\end{equation}

\subsection{Optimal reachability with obstacle avoidance}\label{optimalcont}
We interpret the general static optimal reachability problem, discussed in Section \ref{s2}, in the framework of soft robotics. The control set is $U=[-1,1]$, the configuration is $q(s;u)$ is a solution of the control ordinary differential equation \eqref{reducedstationary}. As in the discrete case, we choose a control quadratic running cost $\ell(u):=u^2$.  Then, given an obstacle $\Omega \subset \RR^2$ and a target point $q^*\in \RR^2\setminus \Omega$, the general problem \eqref{pbgen} reads
\begin{equation}\label{contactfunctionalstationary}
\min \mathcal J,\quad\text{subject to \eqref{reducedstationary} and to $|u|\leq 1$\,,}
\end{equation} 
where
 \begin{equation}
\begin{split}
    \mathcal J(q,u):=&\frac12\int_0^1 u^2(s) ds+\frac{1}{2\delta}|q(1)-q^*|^2+ \frac{1}{2\tau}\int_0^1 
 \mathbf d^2(q(s),\Omega^c)ds,
\end{split}\end{equation}
with $\delta,\tau>0$.
We recall that the three components of the above cost functional respectively represent: a quadratic cost on the controls, a 
tip-target distance, and an integral term vanishing if and only if no interpenetration with the obstacle $\Omega$ occurs. Similarly to the discrete case,  the input-to-state map \eqref{control2state} allows to reduce  $\mathcal J$ to a functional depending on the control $u$ only.  

\subsubsection{Numerical simulations.} 
\begin{figure}[ht!]
	\centering
	\begin{tabular}{c}
		\includegraphics[width=.28\textwidth]{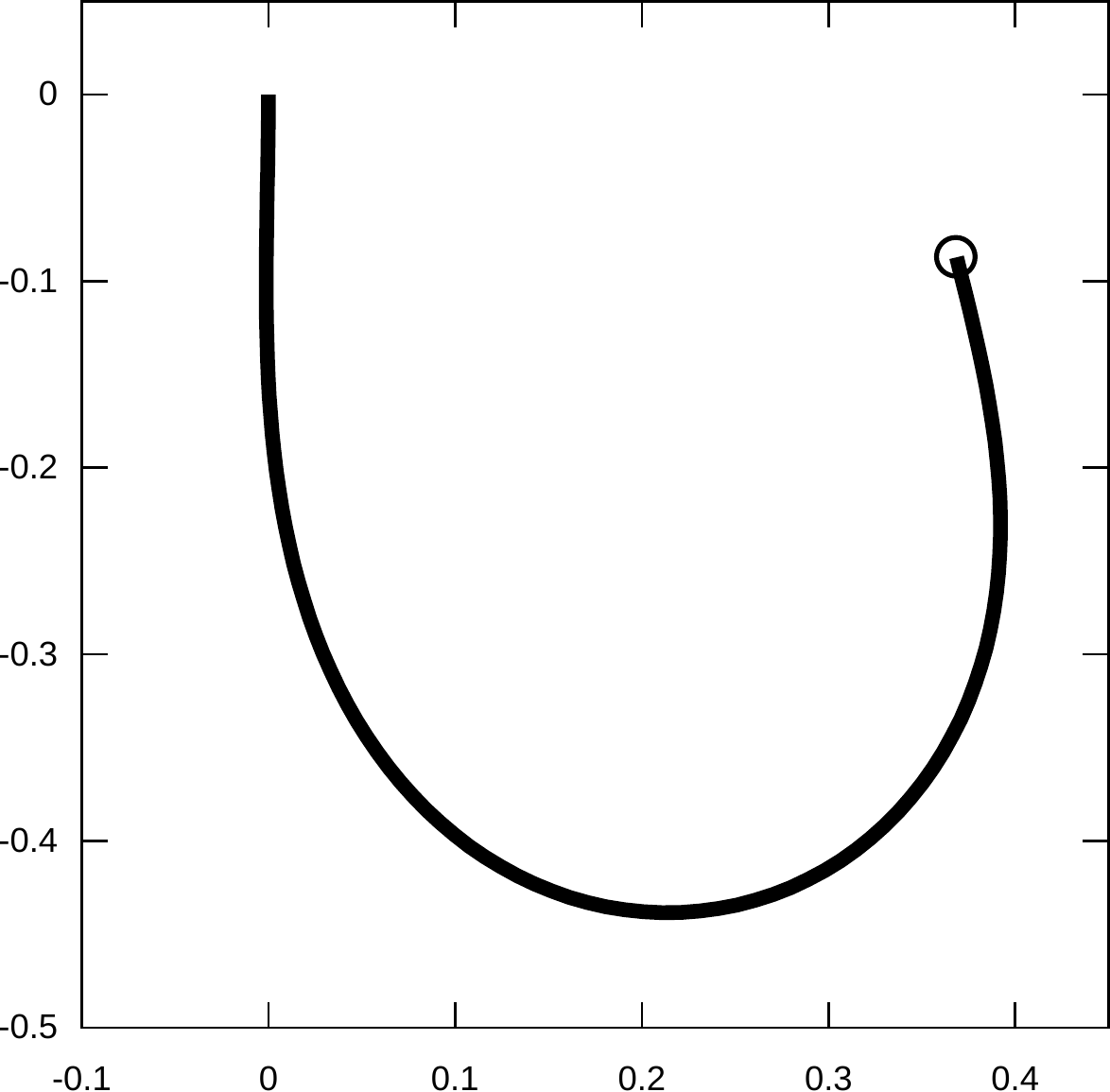}\\
		(a.1)\\
		\includegraphics[width=.28\textwidth]{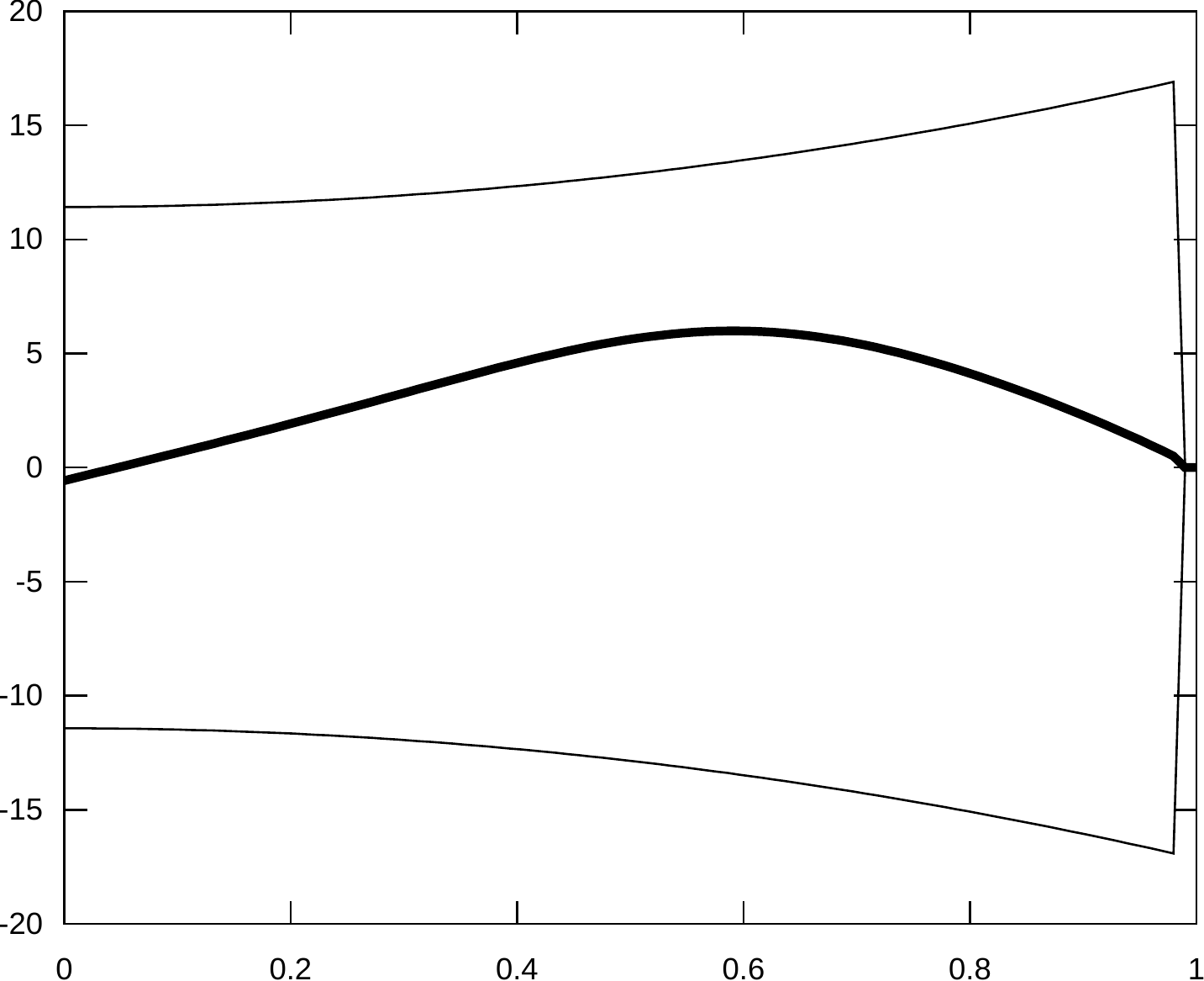}\\
		(b.1)	
	\end{tabular}
	\begin{tabular}{c}
		\includegraphics[width=.28\textwidth]{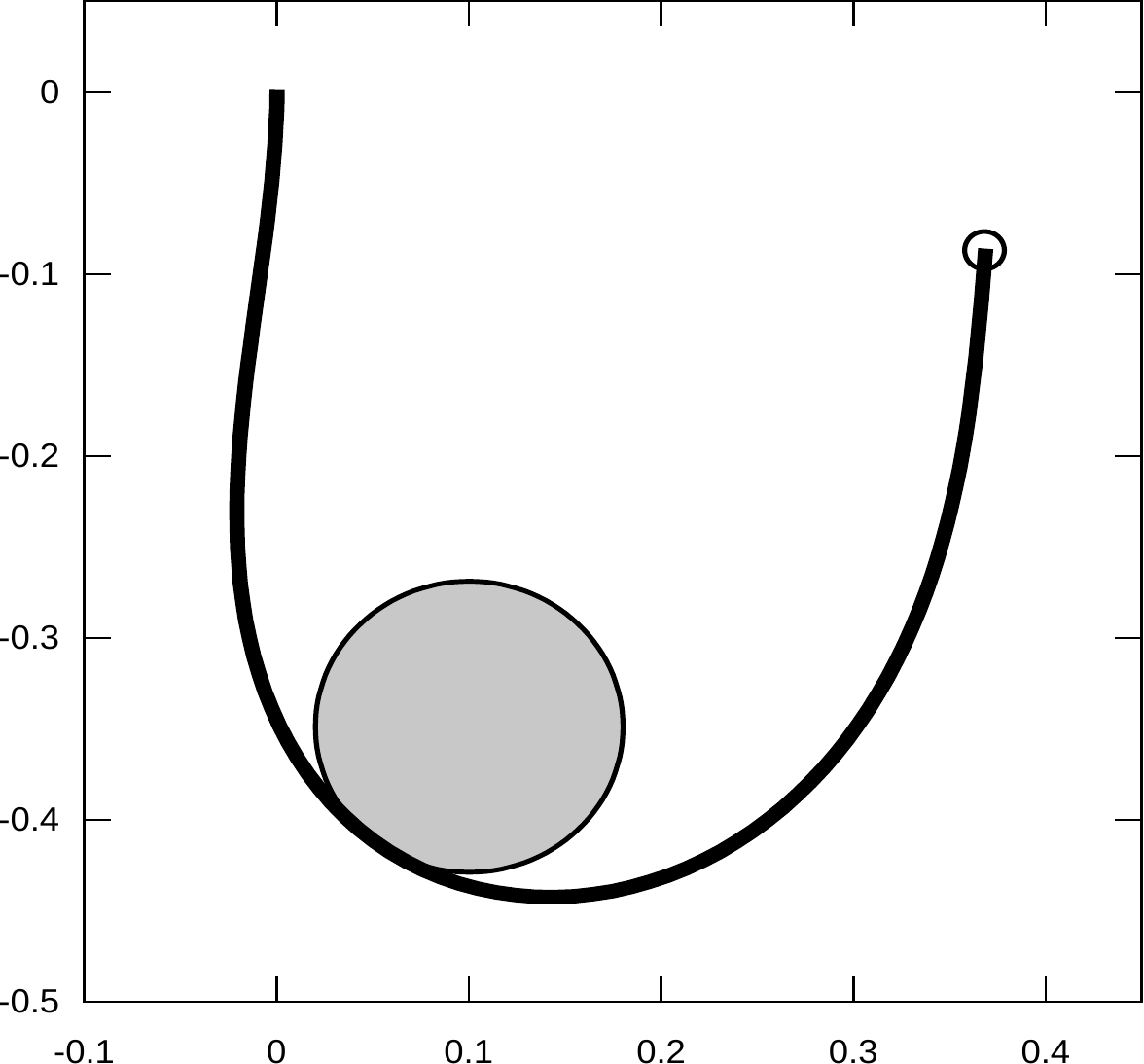}\\
		(a.2)\\
		\includegraphics[width=.28\textwidth]{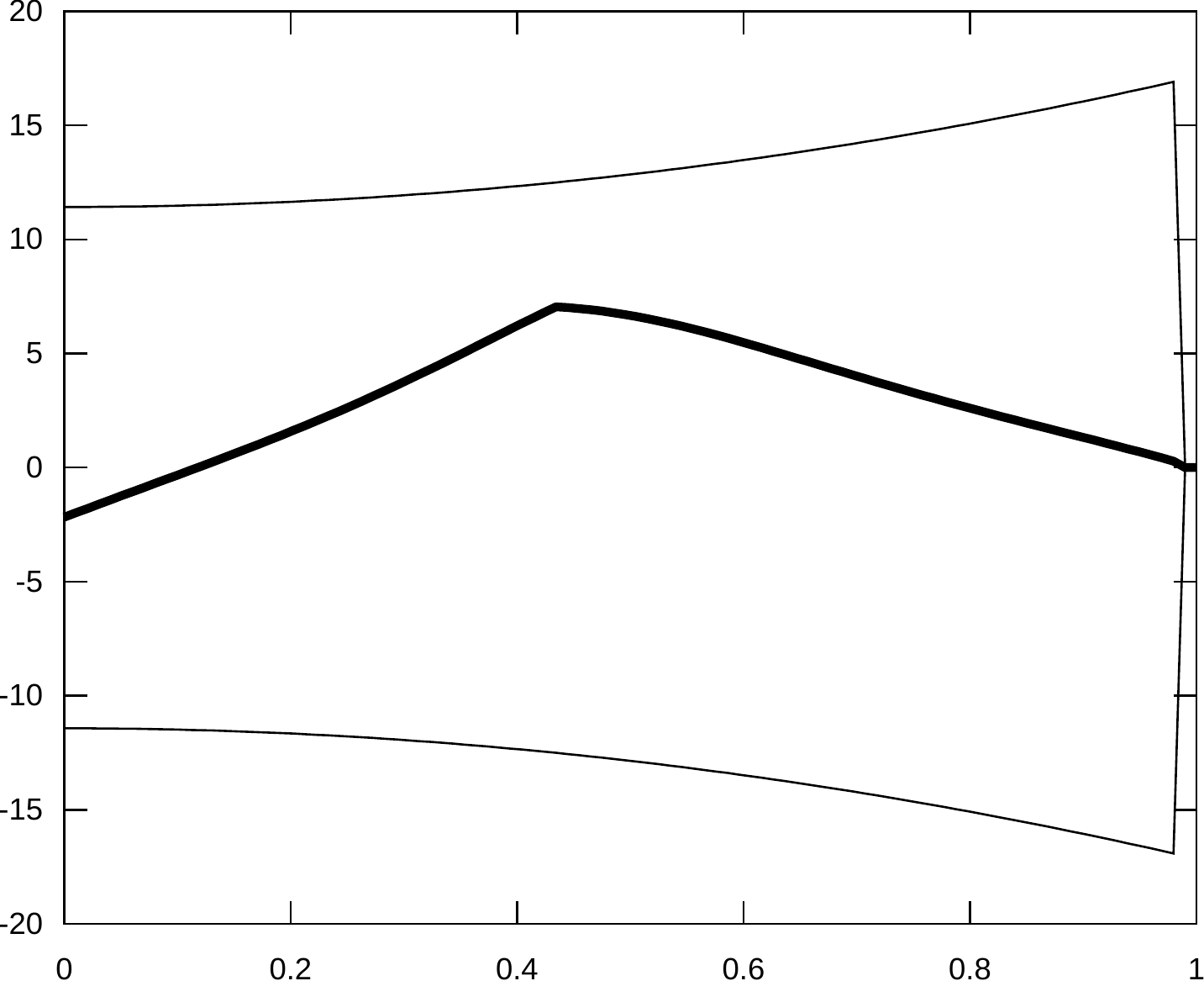}\\
		(b.2)	
	\end{tabular}
		\begin{tabular}{c}
		\includegraphics[width=.28\textwidth]{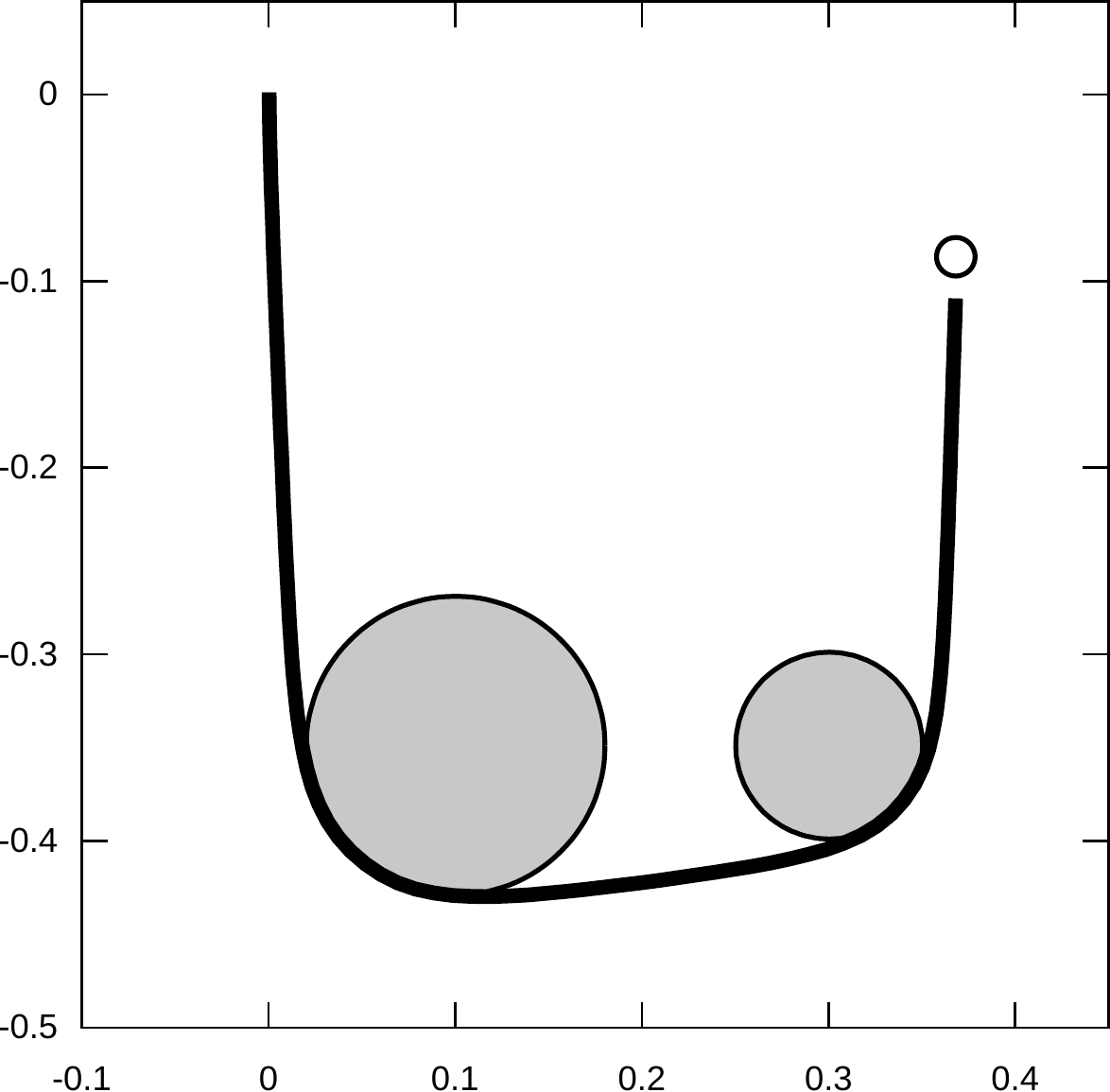}\\
		(a.3)\\
		\includegraphics[width=.28\textwidth]{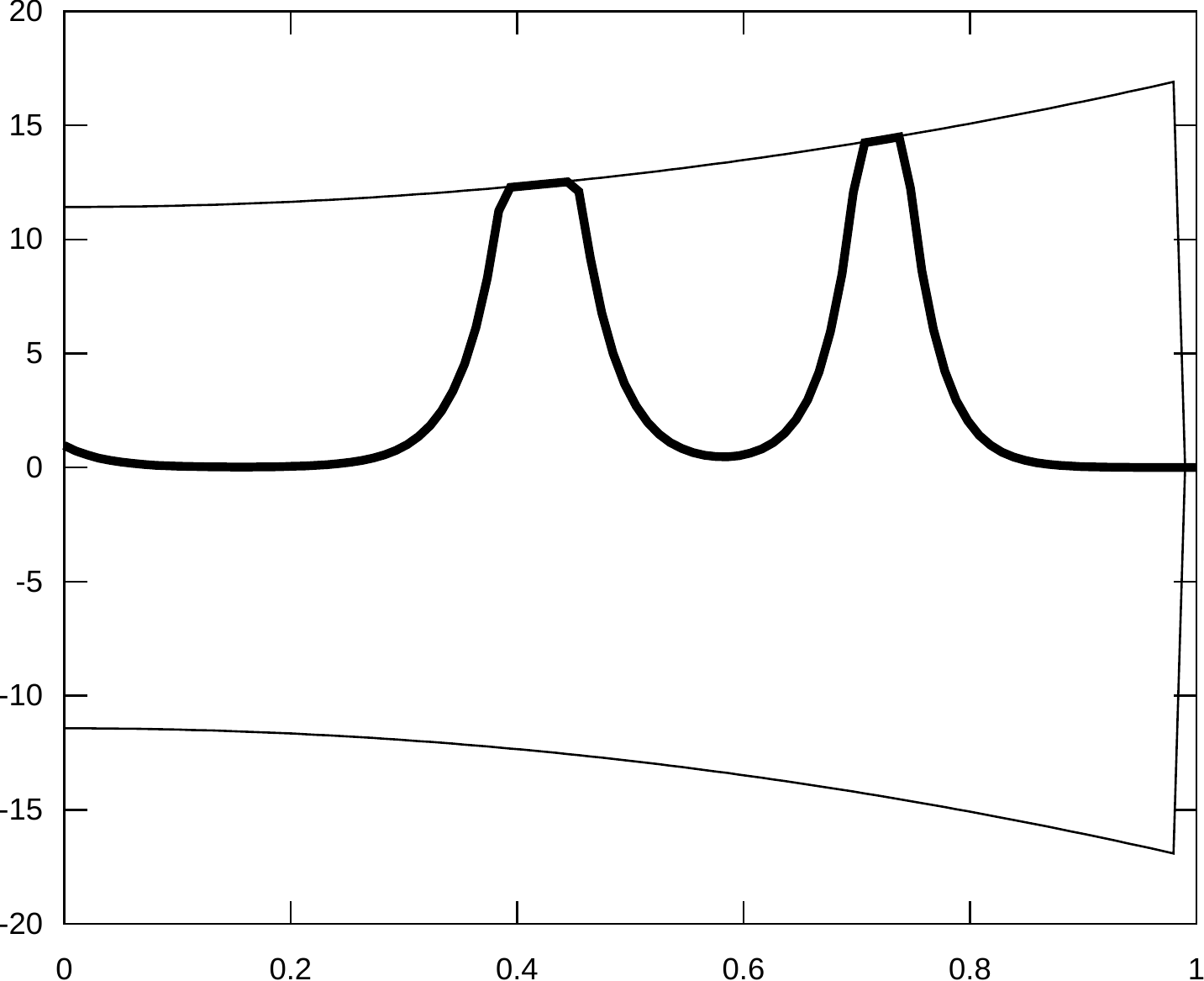}\\
		(b.3)	
	\end{tabular}
	\begin{tabular}{c}
		\includegraphics[width=.28\textwidth]{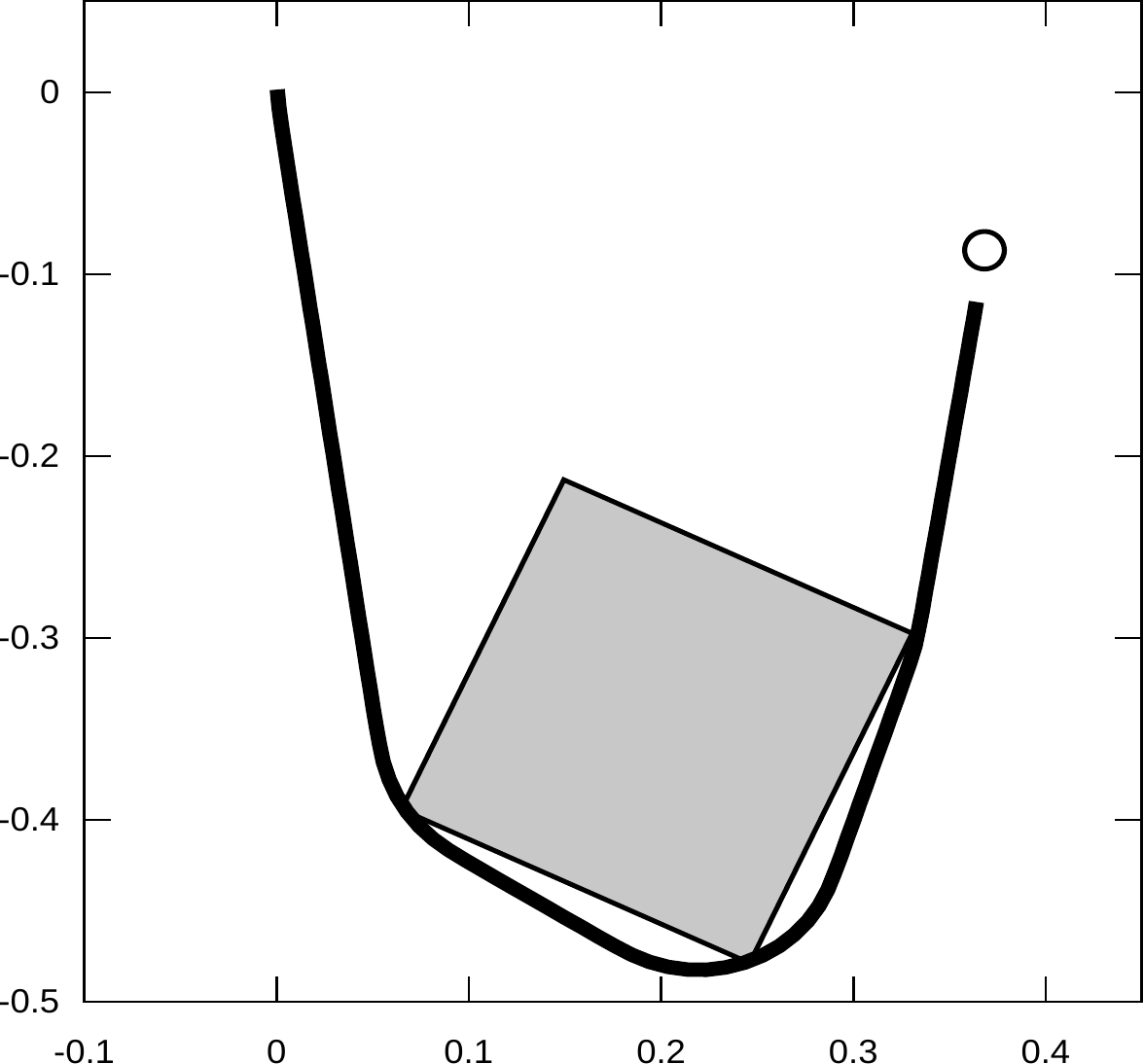}\\
		(a.4)\\
		\includegraphics[width=.28\textwidth]{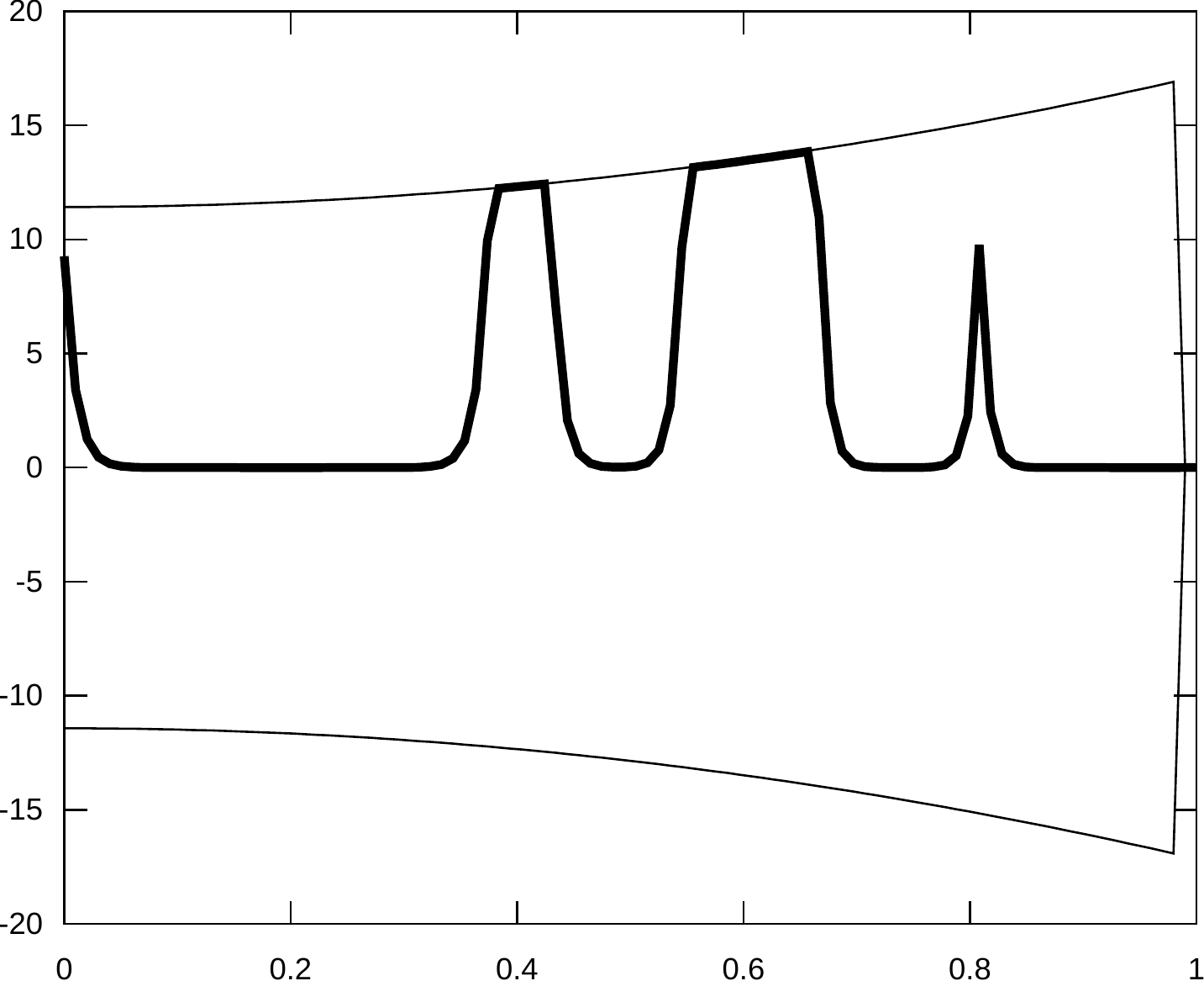}\\
		(b.4)	
	\end{tabular}
	\begin{tabular}{c}
		\includegraphics[width=.28\textwidth]{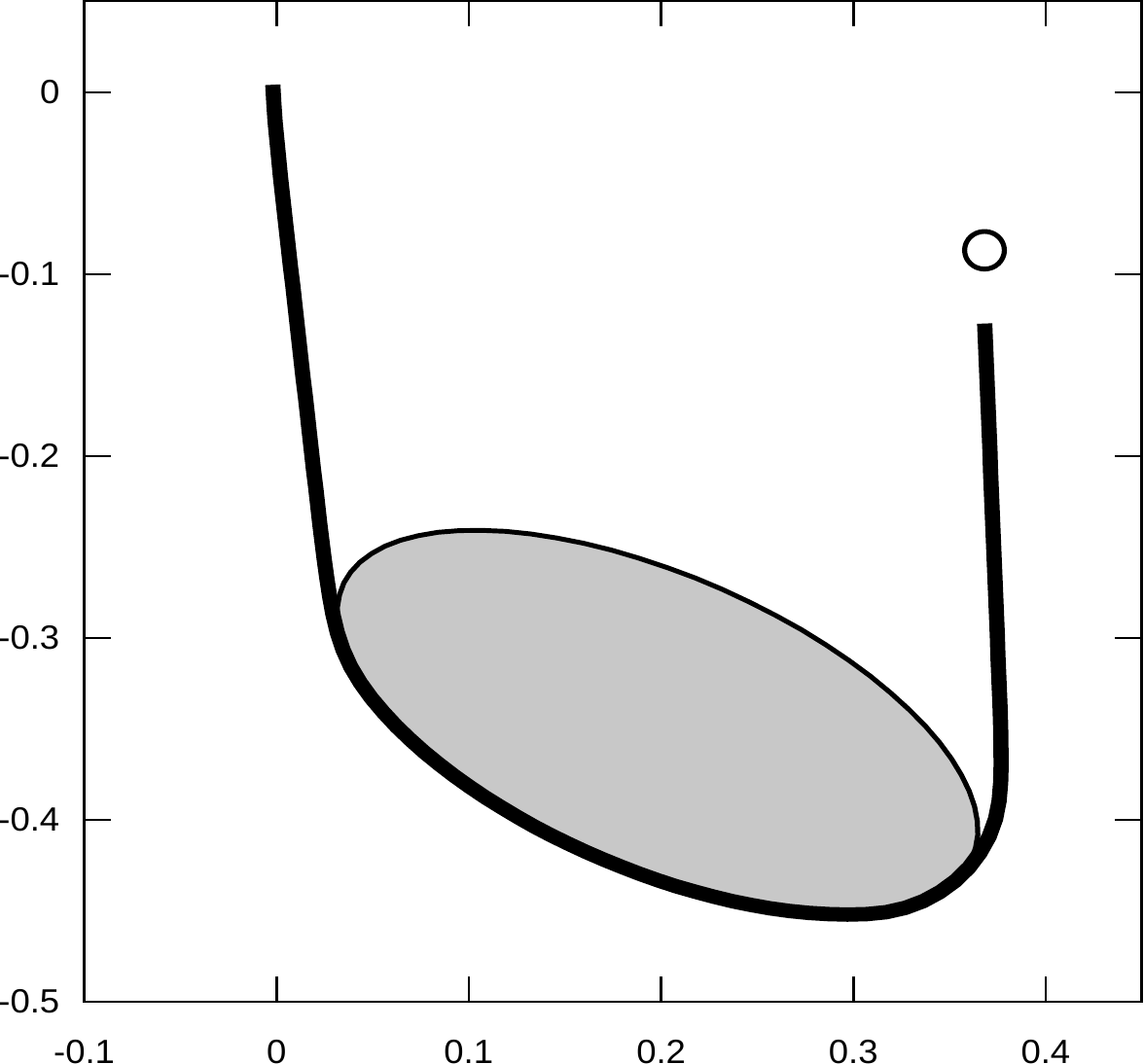}\\
		(a.5)\\
		\includegraphics[width=.28\textwidth]{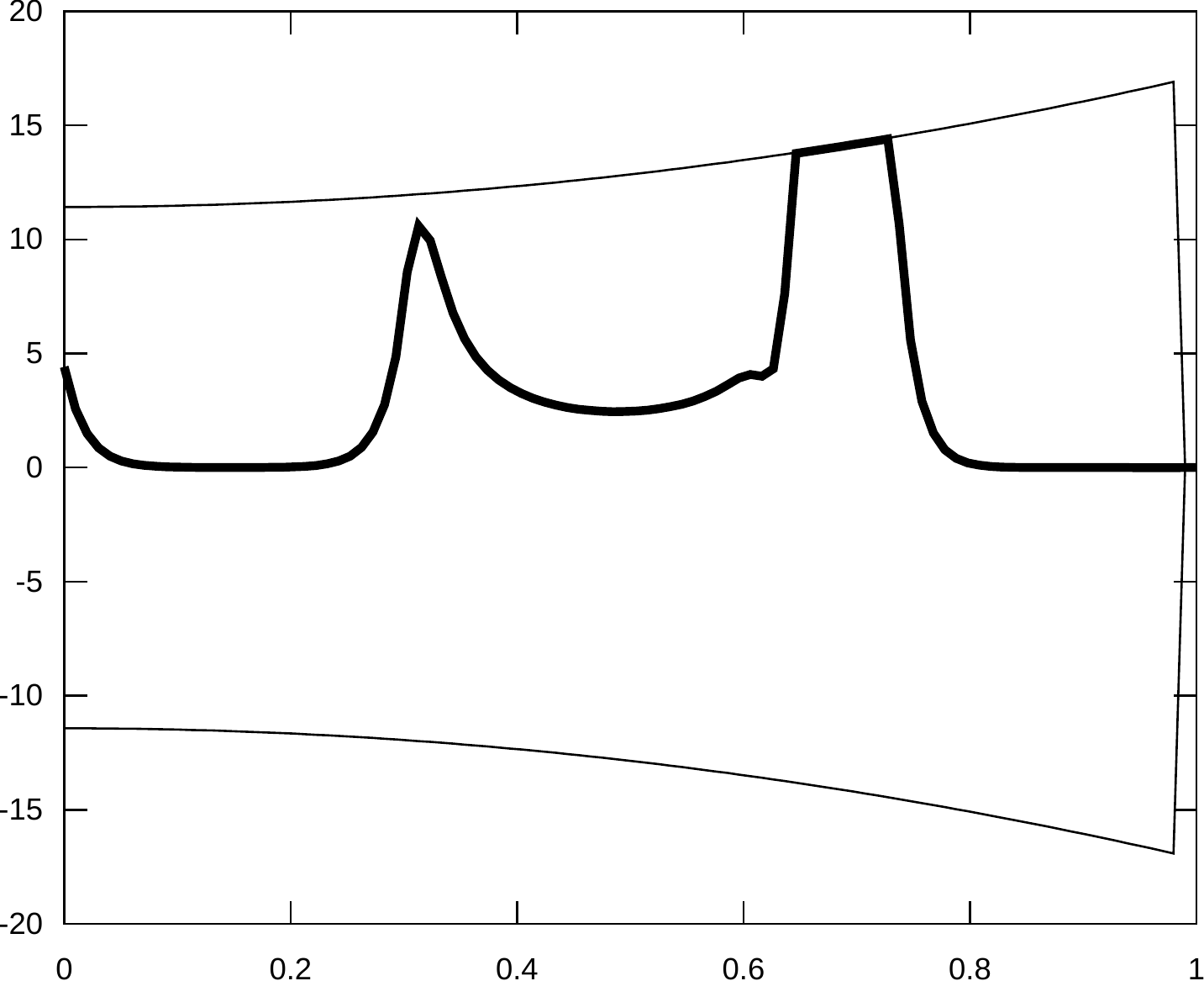}\\
		(b.5)	
	\end{tabular}
		\begin{tabular}{c}
		\includegraphics[width=.28\textwidth]{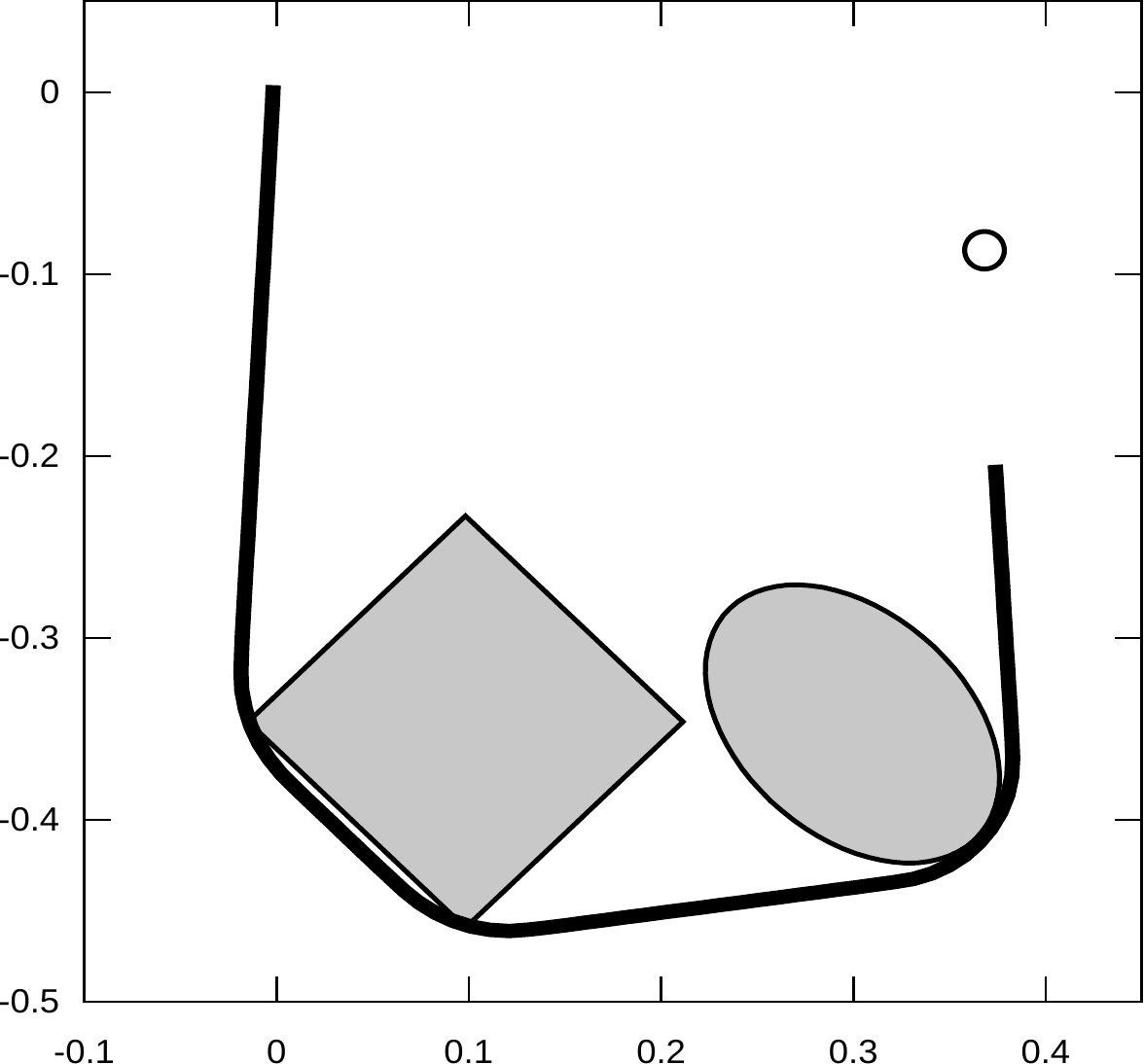}\\
		(a.6)\\
		\includegraphics[width=.28\textwidth]{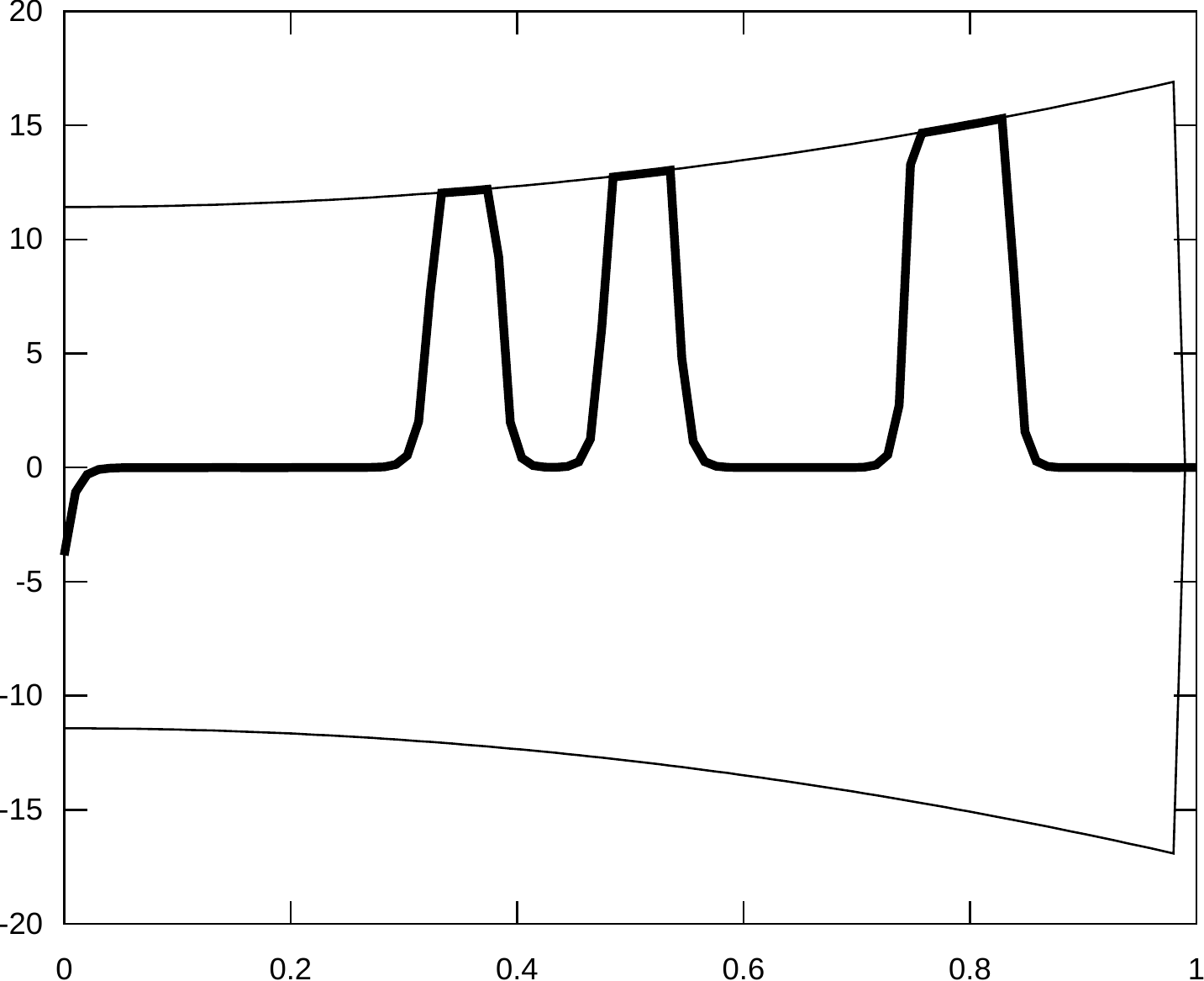}\\
		(b.6)	
	\end{tabular}
	\caption{\label{stationary3}  Figures labeled with (a.$\sharp$) represent the solution $q$ of Test $\sharp$, in (b.$\sharp$) the related signed curvature $\kappa(s)$ (bold line) and curvature constraints $\pm\bar \omega$ (thin lines).} 
\end{figure}

Discretization and optimization are performed as in the case of hyper-redundant manipulators, using  quadrature rules  to approximate the integrals appearing in the input-to-state map \eqref{control2state} and in the functional \eqref{contactfunctional}. For the sake of comparison, we adopt the same obstacle settings of the discrete case, reported in Table \ref{controltable}. The other global parameter settings are in Table \ref{parametertable}. We note that in Test 1 and Test 2 the target is reached and the optimal controlled curvature $\kappa$ is far below the fixed threshold $\bar\omega$ -- see Figure \ref{stationary3}(a.1-2) and Figure \ref{stationary3}(b.1-.2). The remaining tests displayed in Figure \ref{stationary3} show more clearly the impact of curvature and obstacle avoidance constraints on the optimization process: the optimal configuration fails in reaching the target. \begin{table}[t]
  \begin{center}	\begin{tabular}{|l|l|}\hline
  	
  		Parameter description&Setting\\
  		\hline Quadrature nodes &$N=100$\\
  		Discretization step &$\Delta_s=1/N=0.01$\\
  		\hline 
  		bending moment &$\varepsilon(s)=10^{-1}(1-0.9s)$\\
  		curvature control&$\mu(s)=1-0.9s$\\ 
  		penalty&\\
  		curvature constraint&$\omega(s)=2\pi(2+s^2)$\\
  		\hline
  		target point&  $q^*=(0.368,-0.085)$\\
  		target penalty &$\delta=10^{-8}$\\
  		obstacle penalty &$\tau=10^{-10}$\\
  		\hline
  	\end{tabular}\end{center}
  	\caption{Global parameter settings for the soft manipulator.  \label{parametertable}}
  \end{table}

 \bibliographystyle{splncs04}
 \bibliography{vurpo}
\end{document}